\documentclass{amsart}
\usepackage{amssymb,mathrsfs,amsmath,dsfont}
\usepackage{tikz}
\usetikzlibrary{arrows,positioning}
\usepackage{capt-of}% It helps you to use caption for tikzpictures as well.
\usepackage{mathbbol}     % Using these you can use mathbb fonts for lowercase characters as well.
\usepackage{color}
\usepackage{verbatim} %Allows use of \begin{comment} ... \end{comment}
\usepackage[all,cmtip]{xy}
\usepackage{amsbsy}% this helps us to use bold fonts for mathematical symbols as well. We then use $\pmb{any symbol}.$
\usepackage{amsmath}% this allows \tag{any symbol} to be used  for the equations.
\usepackage{nccmath}% this allows you to have equation with position on the left. \begin{fleqn}\begin{flalign}\end{fleqn}\end{fleqn}

\usepackage{hyperref}

\newtheorem{Thm}{Theorem}[section]
\newtheorem{Pro}[Thm]{Proposition}
\newtheorem{cor}[Thm]{Corollary}
\newtheorem{lem}[Thm]{Lemma}
\newtheorem{exa}[Thm]{Example}

\newtheorem{deft}[Thm]{Definition}
\newtheorem{rem}[Thm]{Remark}

\numberwithin{equation}{section}
\newenvironment{prf}{\begin{proof}}{\end{proof}}

\newcommand{\bT}{{\mathbb T}}
\newcommand{\one}{\mathbf{1}}

\newcommand{\cH}{\mathcal{H}}
\newcommand{\cC}{\mathcal{C}}
\newcommand{\cO}{\mathcal{O}}
\newcommand{\cG}{\mathcal{G}}

\newcommand{\cN}{\mathcal{N}}

\newcommand{\R}{{\mathbb R}}
\newcommand{\Z}{{\mathbb Z}}
\newcommand{\C}{{\mathbb C}}

\newcommand{\oline}{\overline}

\newcommand{\res}{\vert}
\newcommand{\Aut}{\mathop{{\rm Aut}}\nolimits}
\newcommand{\conv}{\mathop{{\rm conv}}\nolimits}
\newcommand{\PU}{\mathop{{\rm PU}}\nolimits}
\newcommand{\Herm}{\mathop{{\rm Herm}}\nolimits}
\newcommand{\Ind}{\mathop{{\rm Ind}}\nolimits}
\newcommand{\OO}{\mathop{\rm O{}}\nolimits}
\newcommand{\Pin}{\mathop{{\rm Pin}}\nolimits}
\newcommand{\Spin}{\mathop{{\rm Spin}}\nolimits}
\newcommand{\SO}{\mathop{{\rm SO}}\nolimits}
\renewcommand{\:}{\colon}
\newcommand{\cent}{\mathfrak{cent}}
\newcommand{\End}{\mathop{{\rm End}}\nolimits}
\newcommand{\Sym}{\mathop{{\rm Sym}}\nolimits}
\newcommand{\der}{\mathop{{\rm der}}\nolimits}
\newcommand{\diag}{\mathop{{\rm diag}}\nolimits}
\newcommand{\id}{\mathop{\rm id}}

\newcommand{\urad}{\mathop{{\rm urad}}\nolimits}
\newcommand{\subeq}{\subseteq}

\newcommand{\eps}{\varepsilon}
\renewcommand{\hat}{\widehat}
\renewcommand{\tilde}{\widetilde}

\def\a1s{a_1,\cdots, a_s}
\def\a{\alpha}

\def\andd{\quad\hbox{and}\quad}

\def\ad{\hbox{ad}}

\def\b{\beta}

\def\bl4{B_{\ell\geq4}}
\def\sc{\mathscr{C}}
\def\cc{{\mathcal C}}
\def\fc{\mathfrak{c}}
\def\bbbc{{\mathbb C}}

\def\d{\delta}

\def\fe{\mathfrak{E}}

\def\bbbf{\mathbb{F}}

\def\psu{\mathfrak{psu}}
\def\pq{\mathfrak{pq}}
\def\fg{\mathfrak{g}}
\def\g{\mathfrak{g}}

\def\hH{\mathscr{H}}
\def\hN{\mathscr{N}}
\def\hM{\mathscr{M}}

\def\hK{\mathscr{K}}
\def\hB{\mathscr{B}}

\def\fh{\mathfrak{h}}
\def\su{\mathfrak{su}}

\def\fk{\mathfrak{k}}

\def\lam{\lambda}
\def\Lam{\Lambda}

\def\ep{\epsilon}
\def\fm{(\cdot,\cdot)}

\def\fn{\frak{n}}

\def\bbbr{{\mathbb R}}

\def\1k{\frac{1}{k}}
\def\op{\oplus}
\def\ot{\otimes}

\def\la{\langle}
\def\ra{\rangle}

\def\sub{\subseteq}
\def\sg{\sigma}

\def\i{{\mathcal I}}

\def\fu{\frak{u}}
\def\fsu{\frak{su}}
\def\fpsu{\frak{psu}}

\def\bbbz{{\mathbb Z}}

\def\1il{1\leq i\leq\ell}

\renewcommand{\hat}{\widehat}

\begin{document}
\centerline{\bf }
\vspace{0.7cm}

\title{Current superalgebras and unitary representations}
%\thanks{Address: Department of Mathematics, University of Isfahan, Isfahan, Iran,
%P.O.Box 81745-163, and School of Mathematics, Institute for Research in
%Fundamental
%Sciences (IPM), P.O. Box: 19395-5746, Tehran, Iran.\\
%Email address: ma.yousofzadeh@sci.ui.ac.ir \& ma.yousofzadeh@ipm.ir.\\
%}
\maketitle
\centerline{Karl-Hermann Neeb\footnote{Department Mathematik, FAU, Cauerstrasse 11, 91058 Erlangeen,
Germany. Email address: neeb@math.fau.de.} and Malihe Yousofzadeh\footnote{Department of Mathematics, University of Isfahan, Isfahan, Iran,
P.O.Box 81745-163, and School of Mathematics, Institute for Research in
Fundamental
Sciences (IPM), P.O. Box: 19395-5746, Tehran, Iran.
Email address: ma.yousofzadeh@sci.ui.ac.ir \& ma.yousofzadeh@ipm.ir.}}
\vspace{2cm}

%\today

\begin{abstract}
{In this paper we determine the projective unitary representations of finite
dimensional Lie supergroups whose underlying Lie superalgebra is $\g = A \otimes \fk$, where
$\fk$ is a compact simple Lie superalgebra and
$A$ is a supercommutative associative (super)algebra; the crucial case is
when $A = \Lambda_s(\R)$ is a Gra\ss{}mann algebra.
Since we are interested in projective representations,
the first step consists in determining the cocycles defining the corresponding
central extensions. Our second main result asserts that,
if $\fk$ is a simple compact Lie superalgebra with $\fk_1\neq \{0\}$,
then each (projective) unitary representation of $\Lam_s(\bbbr)\ot \fk$
factors through a (projective) unitary representation of $\fk$
itself, and these are known by Jakobsen's classification.
If $\fk_1 = \{0\}$, then we likewise reduce the classification
problem to semidirect products of compact Lie groups $K$ with a
Clifford--Lie supergroup which has been studied by Carmeli, Cassinelli, Toigo
and Varadarajan.
}\end{abstract}

\tableofcontents

\section{Introduction}

In a similar fashion as projective unitary representations
$\pi \: G \to \PU(\hH)$ of a Lie group $G$ implement symmetries
of quantum systems modelled on a Hilbert space $\hH$,
projective unitary representations of Lie supergroups
implement symmetries of super-symmetric quantum systems \cite{CCTV}.
Here the Hilbert space is replaed by a super Hilbert space
$\hH = \hH_0 \oplus \hH_1$, i.e., a direct sum of two subspaces
corresponding to a $\Z_2$-grading of $\hH$.
We deal with Lie supergroups as Harish--Chandra pairs
$\cG = (G,\g)$, where $G$ is a Lie group and
$\g = \g_0 \oplus \g_1$ is a Lie superalgebra, where
$\g_0$ is the Lie algebra of $G$, and we have an {\it adjoint action}
of $G$ on $\g$ by automorphisms of the Lie superalgebra $\g$ extending
the adjoint action of $G$ on its Lie algebra $\g_0$.

The concept of a unitary representation of a
Lie supergroup $\cG$ consists of a
unitary representation $\pi$ of $G$ by grading preserving
unitary operators and  a representation $\chi_\pi$
of the Lie superalgebra $\g$ on the dense subspace
$\hH^\infty$ of smooth vectors for $G$ such that natural
compatibility conditions are satisfied
(see Definition~\ref{def:unirep-super} for details).
To accomodate the fact that the primary interest lies in projective
unitary representations, one observes that projective
representations lift to unitary representations of central extensions
by the circle group $\bT$ acting on $\hH$ by scalar multiplication.
Having this in mind, one can directly study unitary representations
of central extensions (see \cite{JN} for more details on this passage).

The corresponding classification problem splits into two layers.
One is to determine the even central extensions of a given
Lie supergroup $\cG$ and the second consists of determining for these
central extensions the corresponding unitary representations.

The existence of an invariant measure implies that for any finite
dimensional Lie group $G$, unitary representations exist in abundance,
in particular the natural representation on $L^2(G)$ is injective.
This is drastically different for Lie supergroups, for which
all unitary representations may be trivial. The reason for this is that,
for every unitary representation $\chi \: \g \to \End(\hH^\infty)$
and every odd element $x_1 \in \g_1$, the  operator
$-i\chi([x_1, x_1])$ is non-negative. This imposes serious positivity
restrictions on the representations on the even part $\g_0$, namely that
$-i\chi(x) \geq 0$ for all elements in the closed convex cone
$\sc(\g) \subeq \g_0$ generated by all brackets $[x_1,x_1]$ of odd elements.
Accordingly, $\g$ has no faithful unitary representation if the cone
$\sc(\g)$ is  not pointed (cf.\ \cite{NS}). Put differently, the kernel of
any unitary representation contains the ideal
$\urad(\g)$ of the Lie superalgebra $\g$ generated by the linear subspace
$\fe := \sc(\g) \cap -\sc(\g)$ of $\g_0$ and all those elements
$x \in \g_1$ with $[x,x] \in \fe$.

A particularly simple but nevertheless important class of Lie superalgebras
are the {\it Clifford--Lie superalgebras} $\g$ for which
$[\g_0, \g] = \{0\}$ ($\g_0$ is central), so that the Lie bracket of $\g$ is
determined by a symmetric bilinear map $\mu \: \g_1 \times \g_1 \to \g_0$.
If $\g_0 = \R$ and the symmetric bilinear form $\mu$ is indefinite,
then $\g$ has no non-zero unitary representations.

In \cite{AN} the authors have determined the structure of
finite dimensional Lie superalgebras $\g$ for which finite dimensional
unitary representations exist. This property implies in particular
that $\g$ is {\it compact} in the sense that $e^{\ad \g_0} \subeq
\Aut(\g)$ is a compact subgroup, but, unlike the purely even case,
this condition is not sufficient for the existence of
finite dimensional projective unitary
representations. In particular, it is shown in \cite{AN} that only four
families of simple compact Lie superalgebras have
finite dimensional projective unitary representations:
$\su(n|m;\C), n \not=m$, $\psu(n|n;\C)$, $\fc(n)$ and $\pq(n)$
(see Subsection~\ref{subsec:5.2} for details).

In this paper we take the next step by considering
{\it current Lie superalgebras} $\g = A \otimes \fk$, where
$\fk$ is a compact simple Lie superalgebra and
$A$ is a supercommutative associative (super)algebra
and study the projective unitary representations, resp., the unitary representations of
central extensions of these Lie superalgebras.
Since we are interested in the phenomena caused by the superstructure,
the main interest lies in algebras $A$ generated by their odd part $A_1$.
As the supercommutativity implies that the squares of odd elements in $A$
vanish, any such $A$ is a quotient of a Gra\ss{}mann algebra.
Therefore the main point is to understand current {super}algebras of the form
$\Lambda_s(\R) \otimes \fk$, where $\fk$ is a compact simple Lie superalgebra
and $\Lambda_s(\R)$ is the Gra\ss{}mann algebra with $s$ generators.

Our main result are the following. As explained below, we first have to
understand the structure of the central extensions, resp., of the
even $2$-cocycles. This is described in
Section~\ref{sec:4} and works as follows.
Suppose that $\kappa$ is a non-degenerate invariant symmetric bilinear form
on $\fk$ which is invariant under all derivations of $\fk$. Then
any $D \in \der(\fk)$ and any linear map $f \: A \to \R$ leads to a $2$-cocycle
\[ \eta_{f,D}(a \otimes x, b \otimes y) :=
(-1)^{|b||x|}f(ab)\kappa(Dx,y).\]
There is a second class of natural cocycles on $A \otimes \fk$. To describe it,
we call a bilinear map
$F\colon A\times A\longrightarrow \R$ a {\it Hochschild map} if
\[ F(a, b)=-(-1)^{|a||b|}F(b ,  a) \quad \mbox{ and }  \quad
F(ab ,  c)=F(a, bc)+(-1)^{|b||a|}F(b ,  ac) \]
hold for  $a,b,c\in A$. If $S \: \fk \to \fk$ is $\kappa$-symmetric
and contained in the centroid $\cent(\fk)$, i.e., it commutes with
all right brackets, then
\[ \xi_{F,S}(a \otimes x, b \otimes y) :=
(-1)^{|b||x|}F(a, b)\kappa(Sx,y) \]
also defines a $2$-cocycle.
Our first main result Theorem~\ref{cor1} asserts that each \break
{$2$-cocycle} on $A \otimes \fk$
is equivalent to a finite sum of cocycles of the form
$\eta_{f,D}$ and~$\xi_{F,s}$.

Our second main result is Theorem~\ref{unitary-rep}
which asserts that, if $\fk$ is a simple compact Lie superalgebra with $\fk_1\neq \{0\}$,
then each unitary
representation of $\fg=\Lam_s(\bbbr)\ot \fk$
factors through the quotient map
$\eps \otimes \id_\fk \:  \g \to  \fk$ corresponding to the augmentation
homomorphism $\eps \: \Lambda_s(\R) \to \R$.
This result shows that, if $\fk$ is not purely even, the passage to the
current {super}algebra does not lead to more unitary representations
than what we have seen in \cite{AN} for simple compact Lie superalgebras.
Note that their irreducible representations have been determined
in terms of highest weights by Jakobsen \cite{Jak}.
For an argument that all irreducible unitary representations are of this
form, see \cite{NS}.

This leaves us with the case where $\fk = \fk_0$ is a (purely even) compact Lie
algebra. If $\Lambda_s^+(\R) = \ker \eps$, then
\[ \g \cong (\Lambda_s^+(\R) \otimes \fk) \rtimes \fk \]
is a semidirect sum of the compact Lie algebra $\fk$ and the ideal
$\g^+ := \Lambda_s^+(\R) \otimes \fk$. In Theorem~\ref{urad} we show that
every unitary representation of any central extension $\hat\g$ of $\g$
annihilates the ideal $I := \Lambda_s^{> 3}(\R) \otimes \fk$, resp., its central
extension~$I$. As the quotient $\hat\g/ I$ is a semidirect product
$\hat\fn \rtimes \fk$, where $\hat\fn$ is a Clifford--Lie superalgebra,
the classification of the unitary representations of
the corresponding Lie supergroup $\hat\cN \rtimes K$ can be determined
with the methods developed in details in \cite{CCTV}. We provide a
detailed description of these results in Appendix~\ref{sec:3}.
Theorem~\ref{main1} contains the classification of irreducible
unitary representations of any semidirect product supergroup of the form
$\cG = \cN \rtimes K$, where $K$ is a compact Lie group.
As we have seen above, this combined with the other results provides a complete
description of the irreducible unitary representations of current {super}algebras
of the form $A \otimes \fk$, where $\fk$ is a compact simple Lie superalgebra.

The structure of the paper is as follows.
We collect preliminaries in Section~\ref{sec:2}.
In Section~\ref{sec:4} we then turn to the central extensions,
resp., the $2$-cocycles of current {super}algebras $\g = A \otimes \fk$,
culminating in Theorem~\ref{cor1}.
In Section~\ref{sec:5} our  description of the cocycles is used to
determine the unitary representations of central extensions of the
current {super}algebras $A \otimes \fk$. Here we first turn to the case where
$\fk$ is  compact with $\fk_1 = \{0\}$ and then to the case where
$\fk$ is a simple compact Lie superalgebra with $\fk_1 \not=\{0\}$.
In both cases we reduce the classification problem to situations for which
the solutions are known. In the first case we end up with
semidirect products $\fn \rtimes \fk$ covered by \cite{CCTV} and in the
second case with central extensions of $\fk$ itself. Throughout the paper, we always assume that Lie superalgebras under consideration  are finite dimensional.
%We conclude Section~\ref{sec:5} with the observation
%that for a compact simple Lie superalgebra $\fk$,
%the tangent algebras $\g = \Lambda_1(\R) \otimes \fk = \R[\eps] \otimes \fk$
%may have central extensions for which the cone $\sc(\hat\g)$ is pointed,
%although all unitary representations annihilate the ideal generated by
%$\eps \otimes \fk$.

\section{Preliminaries}
\label{sec:2}

In this section we provide precise definitions required for
unitary representations of Lie supergroups.

A pre-Hilbert space  $(\hH,\la\cdot,\cdot\ra)$ is called  a {\it pre-Hilbert superspace} if $\hH=\hH_0\op\hH_1$ is a superspace such that $\la\hH_0,\hH_1\ra=\{0\}.$ A pre-Hilbert superspace $(\hH,\la\cdot,\cdot\ra)$ is called a   {\it Hilbert superspace} if $\hH$ is a complete space with respect to the metric induced by $\la\cdot,\cdot\ra$.
For a pre-Hilbert superspace $(\hH,\la\cdot,\cdot\ra),$ an endomorphism $T\in \End(\hH)$ is said to be {\it symmetric with respect to   the inner product $\la\cdot,\cdot\ra$} if $\la Tx,y\ra=\la x,Ty\ra$ for all $x,y\in \hH;$ it is called
{\it nonnegative}, denoted by  $T\geq 0,$ if \break $\la Tx,x\ra\geq0$ for all $x\in\hH.$
A linear homogeneous  isomorphism  between two pre-Hilbert superspaces,
preserving  the corresponding inner products, is called {\it unitary.}
We write ${\rm Aut}(\hH)$ (resp. ${\rm Aut}(\hH)_{\rm even}$) for
the group  of all  unitary (resp.\ even) {\it unitary automorphisms} of~$\hH$.

\begin{deft} {\rm
(i) Let $(\hH,\la \cdot,\cdot\ra)$ be  a Hilbert superspace and $G$ be a {finite dimensional Lie group}. A {\it unitary representation} of $G$ on $\hH$
is a pair $(\pi, \hH)$, where
$\pi:G\longrightarrow {\rm Aut}(\hH)$ is a group homomorphism such that,
for each $\zeta\in\hH,$ the orbit map
\begin{align*}
\pi^\zeta:G\longrightarrow \hH,&\quad
g\mapsto \pi(g)\zeta
\end{align*}
is  continuous.
%Then we call $\hH$ is a {\it unitary $G$-module}.
%A closed subspace $\mathscr{K}$ of $\hH$ is called a {\it unitary submodule} of $\hH$
%if each $\pi(g)$ preserves~$\mathscr{K}.$
%The representation $\pi$ (and correspondingly the module $\hH$) is called {\it irreducible} if $\hH$ does not have nontrivial  unitary submodules.
An element $\zeta\in\hH$ is called a {\it smooth vector} if {$\pi^\zeta$} is a smooth  function
\cite[\S III.3]{K}.  We denote by $\hH^\infty$ the set of all smooth vectors of $(\pi,\hH)$ and recall that it is a dense subset of $\hH$ as $G$ is finite dimensional \cite{G}. As $G$ acts by homogeneous operators, $\hH^\infty$ is a sub-superspace of $\hH.$

(ii) A {\it unitary representation} of a real  Lie superalgebra
$\fg$ in  a   pre-Hilbert superspace  $(\hH,\la\cdot,\cdot\ra)$ is a
real  Lie superalgebra homomorphism
$\chi:\fg\longrightarrow {\rm End}(\hH)$ satisfying
\[ \la \chi(X)(u),v\ra=\la u,-i^{|X|}\chi(X)(v)\ra\;\;\;\;\;\;\;\;\;{\hbox{ for $X\in \fg,\; u,v\in \hH.$}}\]
{We then}  refer to $\hH$ as a {\it unitary $\fg$-module.}

(iii) Two unitary  representations are said to be {\it equivalent}
(or {\it isomorphic}) if their actions  intertwine with  an even unitary operator.}
\end{deft}

 \begin{deft} \label{def:lie-subgrp}
{\rm
(i) A {\it Lie supergroup} is a pair $\cG := (G,\fg),$ in which $G$ is a (finite dimensional) Lie group and $\fg=\fg_0\op\fg_1$ is a  {finite dimensional} real Lie superalgebra such that
\begin{itemize}
\item $\fg_0$ is the Lie algebra of $G,$
\item there is a {smooth} action  $\sg:G\longrightarrow \Aut(\fg)$ of $G$ on $\fg,$
\item the differential of $\sg$ is the adjoint action of $\fg_0$ on $\fg.$
\end{itemize}
We denote the Lie supergroup $(G,\fg)$ by  $(G,\fg,\sg)$  {if we want to emphasise~$\sg.$}

 (ii) A {\it Lie subsupergroup} of a
Lie supergroup $\cG = (G,\fg)$ is a pair $\cH = (H,\fh)$
in which   $H$ is a Lie subgroup  of $G,$   $\fh=\fh_0\op\fh_1$ is a Lie
sub-superalgebra  of $\fg$  and the action of $H$ on $\fh$ is the restriction of the action of $G$ on  $\fg.$ A Lie subsupergroup is called {\it special} if
$\fg_1 = \fh_1$.

(iii) If $(G,\fg)$ is a Lie supergroup, then an {\it automorphism of $(G,\fg)$}
is a pair $(\gamma, \beta) \in \Aut(G) \times \Aut(\fg)$ such that
$\beta\res_{\fg_0}$ coincides with the differential $d\gamma$
and $\beta \sigma(g) \beta^{-1} = \sigma(\gamma(g))$ for $g \in G$.

(iv) If $\cN = (N,\fn)$ is a Lie supergroup, $K$ a Lie group and
\[\alpha  = (\alpha_N, \alpha_\fn)\: K \to \Aut(N,\fn) \]
is a smooth group homomorphism, then
we can form the {\it semidirect product Lie supergroup}
$\cN \rtimes_\alpha K := (N \rtimes_{\alpha_N} K, \fn \rtimes_{\beta_\fn} \fk)$, where
$\beta_\fn \: \fk \to \der(\fn)$ is the derived action of $\fk$ {(via $\a_\fn$)}
by even derivations on the Lie superalgebra $\fk.$ \footnote{ For  a
Lie superalgebra $\fg,$  a linear map $D:\fg\longrightarrow \fg$ is called a {\it derivation} of $\fg$ if, for $x,y\in \fg,$  $D[x,y]=[Dx,y]-(-1)^{|x||y|}[Dy,x].$    The set of derivations of $\fg$ is denoted by $\der(\fg).$}}
\end{deft}

 \begin{deft}\label{def:unirep-super}
{\rm
A {\it unitary pre-representation} of a Lie  supergroup
$(G,\fg,\sg)$ in  a  Hilbert superspace  $(\hH,\la\cdot,\cdot\ra)$ is a
triple $(\pi, \chi_\pi, \hB)$, where
\begin{itemize}
\item $\pi:G\longrightarrow {\rm Aut}(\hH)_{\rm even}$
is a unitary representation {of}  $G,$
\item  $\hB \subseteq \hH$ is a $\pi(G)$-invariant dense subspace contained in the space
$\hH^\infty$  of smooth vectors of $(\pi,\hH),$
\begin{itemize}
\item[\rm(a)] $\chi_\pi:\fg\longrightarrow \End(\hB)$
is a unitary representation of $\fg$ in $\hB,$
\item[\rm(b)]  $\chi_\pi(X)={\rm d}\pi(X)\mid_{\hH^\infty}$ {for $X\in\fg_0,$}
\item[\rm(c)] $\chi_\pi(\sg_g(X))=\pi(g)\chi_\pi(X)\pi(g)^{-1}$ for $X\in \fg_1$
and  $g\in G.$
\end{itemize}
\end{itemize}

A {\it unitary representation} is a unitary pre-representation
for which $\hB = \hH^\infty$ is the full space of smooth vectors of $(\pi,\hH)$.
According to \cite{CCTV} (see also \cite[Lemma~4.4]{MNS}),
every unitary pre-representation
extends uniquely to a unitary representation. {We simply denote a unitary representation $(\pi,\chi,\hH^\infty)$ by $(\pi,\chi).$}

(iv) Suppose that $(G,\fg)$ is a Lie supergroup  and  $(\pi,\chi_\pi)$  is a unitary representation of $(G,\fg)$ in a Hilbert superspace $\hH.$ A closed sub-superspace $\mathscr{K}$ of $\hH$ for which $\pi(g)(\mathscr{K})\sub \mathscr{K}$ and $\chi_\pi(X)(\mathscr{K}^\infty)\sub \mathscr{K}^\infty,$ for all $g\in G$ and $X\in\fg,$  is called
a {\it submodule} of $\hH.$ The unitary representation $\pi$ (and correspondingly the {\it unitary module} $\hH$) is called {\it irreducible} if $\hH$ has no
nontrivial submodule.}
\end{deft}

As we shall need it below,
we recall the construction of unitarily induced representations
in the context of Lie supergroups (see \cite[\S 3]{CCTV} for more details).

\begin{deft}[Induced Representation]
\label{def:ind-rep}
{\rm  Suppose that $\cG = (G,\fg)$ is a Lie supergroup  and $\cH = (H,\fh)$ is a
special Lie subsupergroup of $(G,\fg)$ (Definition~\ref{def:lie-subgrp}).
Suppose $(\rho,\chi_\rho, \hK^\infty)$ is a unitary representation of
$\cH$ and that the (purely even) homogeneous space $H \setminus G$ carries a
$G$-invariant measure $\mu$. Define $\hH$ as
the space of  (equivalence classes of)
measurable functions $f : G \longrightarrow \mathscr{K}$ such that
\begin{itemize}
\item[\rm(a)] for any $g \in G$ and $h \in  H,$  we have $f(hg) = \rho(h)f(g),$
\item[\rm(b)] $\int_{H\setminus G}\|f(g)\|^ 2{d\mu(H g)}< \infty .$
\end{itemize}
In other words, this is the Hilbert space of $L^2$-sections of the
Hilbert bundle $\hK \times_H G$ over $H \setminus G$ associated to the
$H$-principal bundle $G$.

\noindent We define
$\pi:G\longrightarrow {\rm Aut}(\hH)$ by
\[
(\pi(g)f) (g_0) = f(g_0g) \quad\mbox{for } \quad
f\in \hH, g, g_0\in G.\]
Let  $\hB\subeq \hH$ be the subspace of $\hH^\infty$
consisting of all smooth functions \break $f \: G \to \hK$ with
compact support modulo~$H.$ Then $\hB$ is a dense $G$-invariant
subspace of $\hH^\infty$ and  we define
$\chi_{\pi}:\fg\longrightarrow \End(\hB)$ by
\[(\chi_{\pi}(X)f) (g) = \chi_\pi(g\cdot X)f(g) \quad\mbox{ for } \quad
g \in G, X\in \fg_1, f\in \hB. \]
Now  $(\pi,\chi_{\pi})$ defines a unitary
pre-representation of $(G,\fg)$ in $\hH$ and its canonical extension
to a unitary representation is called the {\it induced representation}
and denoted by $(\pi, \chi_\pi) := \Ind_\cH^\cG(\rho, \chi_\rho).$
}\end{deft}

\begin{deft}\label{zero-sq}
{\rm Suppose  $\fg=\fg_0\op\fg_1$ is a finite dimensional
real Lie superalgebra. We write
\[ \sc(\fg) \subeq \fg_0 \]
for the closed convex cone generated by the brackets
$[x,x]$, $ x\in \g_1$ (cf.\ \cite{NS}). The ideal $\urad(\fg)$  of $\fg$
generated by $\fe := \sc(\fg) \cap -\sc(\fg)$
and all those elements $x \in \g_1$ with $[x,x] \in \fe$
is called the {\it unitary radical} of $\fg.$
We say that the convex cone $\sc(\g)$ of $\g$ is {\it pointed}
if $\sc(\g)\cap-\sc(\g)=\{0\}.$ }
\end{deft}

\begin{lem} \label{point-crit} If there exists a linear functional
$\lambda \in \fg_0^*$ with
\[ \lambda([x_1, x_1]) > 0 \quad \mbox{ for } \quad 0 \not= x_1 \in \g_1,\]
then the cone $\sc(\g)$ is pointed.
\end{lem}

\begin{prf} Let $C := \{ x_1 \in\g_1 \: \lambda([x_1, x_1]) = 1\}$.
As  the level set of a positive definite form on $\g_1$, the set $C$ is compact.
Hence $K := \{ [x_1, x_1] \: x_1 \in C\}$ is a compact subset of $\g_0$
contained in the affine hyperplane $\lambda^{-1}(1)$. Therefore the closed
convex cone $\R^+\conv(K) = \sc(\g)$ is pointed.
\end{prf}

The following lemma shows that $\urad(\g)$ is contained in the kernel
of every unitary representation of $\fg$, hence the name.

\begin{lem}\label{necessary}
If $(\chi, \hH)$ is a  unitary representation of the
real Lie superalgebra $\fg$ in a Hilbert superspace, then the following
assertions hold:
\begin{itemize}
\item[\rm(i)] For $x\in \fg_1,$  $\chi(x)=0$ if and only if $\chi([x,x])=0.$
\item[\rm(ii)] $-i\chi(x)\geq 0$ for $x\in \sc(\fg).$
\item[\rm(iii)] $\urad(\fg) \subeq \ker\chi$.
\end{itemize}
\end{lem}

\begin{prf} (i) For $x\in \fg_1$ and $u\in \hH^\infty,$ we have $\chi([x,x])=2\chi(x)^2$ and
\begin{equation}\label{elementary}
-i\la\chi([x,x])u,u\ra=
-2i\la\chi(x)^2u,u\ra=
2\la-i\chi(x)u,-i\chi(x)u\ra\geq0.
\end{equation}

(ii)  follows from (\ref{elementary}).

(iii) First (i) and (ii) imply that
$\sc(\fg)    \cap - \sc(\fg) \subeq \ker \chi$, and as
$\ker\chi$ is an ideal of $\g$, the assertion follows.
\end{prf}

\begin{rem} {\rm(The relation to super-hermitian forms)}\label{rem1}
\\ {\rm(i)  Suppose that $\hH=\hH_0\op\hH_1$ is a superspace equipped with an even  {\it super-hermitian} form $\fm:\hH\times\hH\longrightarrow \bbbc,$ that is,
\begin{itemize}
\item $(\hH_0,\hH_1)=\{0\},$
\item $\fm$ is linear in the first component,
\item $(u,u)>0$ and $-i(v,v)>0$ for $0\neq u\in \hH_0$ and $0\neq v\in\hH_1,$
\item $\overline{(u,v)}=(-1)^{|u||v|}(v,u).$
\end{itemize}
Then  $$\la u,v\ra:=i^{-|u||v|}(u,v)\quad\quad  (u,v\in\hH)$$ defines a
hermitian form for which $\hH_0$ and $\hH_1$ are orthogonal subspaces of $\hH,$
i.e., $(\hH,\la\cdot,\cdot\ra)$ is a Hilbert superspace.
A  homogeneous linear endomorphism \break $T:\hH\longrightarrow \hH$ is called {\it supersvmmetric} with respect to the super-hermitian form $\fm$ if
$$
(Tu,v)=(-1)^{|T||u|}(u,Tv)\quad\quad (u,v\in \hH).
$$
A linear endomorphism $T=T_0+T_1\in \End(\hH)$ is called {\it supersvmmetric} with respect to  $\fm$ if $T_0$ and $T_1$ are supersvmmetric with respect to $\fm.$

The mapping $T\longrightarrow e^{|T|\frac{\pi i}{4}}T$ defines a bijection from the set of supersymmetric linear endomorphisms of $\hH$ with respect to $\fm$ onto symmetric linear endomorphisms of $\hH$ with respect to $\la\cdot,\cdot\ra.$
Moreover, if $\fg$ is a  real Lie superalgebra and $\chi:\fg\longrightarrow \End(\hH)$ is a Lie superalgebra homomorphism, then $\chi$ is  a unitary representation of $\fg$ in the Hilbert superspace $(\hH,\la\cdot,\cdot\ra)$ if and only if, for each $X\in\fg,$ $\chi(X)$ is  a skew-supersymmetric with respect to $\fm,$ i.e.,
\[ (\chi(X)u,v)=-(-1)^{|\chi(X)||u|}(u,\chi(X)v)\quad\quad (u,v\in \hH).\]
}
\end{rem}

\section{Current superalgebras}
\label{sec:4}

In this section, we assume $\bbbf$ is a field of characteristic zero and  unless otherwise mentioned, we consider all vector spaces and tensor products  over $\bbbf.$
The main result of this section is Theorem~\ref{cor1} describing
the structure of the $2$-cocycle of current {super}algebras of the form
$A \otimes \fk$.

\subsection{Invariant forms and $2$-cocycles}
\label{subsec:4.1}

For a Lie superalgebra $\fg$ and a superspace $M,$ a bilinear  map $\omega\colon \fg\times \fg\longrightarrow M$ is called a {\it 2-cocycle with coefficients in $M$} if
\begin{itemize}
\item $\omega(x, y)=-(-1)^{|x||y|}\omega(y, x)$
\item $\omega([x,y], z)=\omega(x,[y,z])-(-1)^{|x||y|}\omega(y,[x,z])$
\end{itemize}
  for all $x,y,z\in\fg.$ The set of all $2$-cocycles with coefficients in $M$ is denoted by $Z^2(\fg,M).$ A $2$-cocycle $\omega$ is called a {\it $2$-coboundary} if there is a linear map $f\colon \fg\longrightarrow M$ with $\omega(x, y)=f([x,y])$ for all $x,y\in \fg.$
   The set of $2$-coboundaries  is denoted by $B^2(\fg,M)$ and the quotient space $$H^2(\fg,M):=Z^2(\fg,M)/B^2(\fg,M)$$ is called the {\it second cohomology} of $\fg$ with coefficients in $M.$
   Two 2-cocycles are called {\it cohomologous} if their {difference}   is a 2-coboundary.

$2$-cocycles of a Lie superalgebra $\fg$  are in correspondence with    its central {extensions:  If} $\fg$ is a Lie superalgebra and $\omega$ is a $2$-cocycle of $\fg$ with coefficients in a superspace $M,$ taking  $\widehat\fg$ to  be the superspace $\fg\op M$ and defining
\begin{align*}
[\cdot,\cdot]_\omega\colon &\widehat\fg\times\widehat\fg\longrightarrow \widehat\fg,\quad\quad
 (x+m,x'+m')\mapsto [x,x']+\omega(x, x')
\end{align*}
for $x,x'\in \fg$ and $m,m'\in M,$ $\widehat\fg$ together with $[\cdot,\cdot]_\omega$ is a Lie superalgebra and the canonical projection map $\pi\colon \widehat\fg\longrightarrow \fg$ is a central extension.
\begin{deft}
{\rm Suppose that $\fg$ is a Lie superalgebra.
A superspace  $M$   together with a bilinear map $\cdot\colon M\times \fg\longrightarrow M$ is called  a {\it right $\fg$-module} if
\[ M_i\cdot \fg_j\sub M_{i+j} \quad \mbox{ and } \quad
a\cdot [x,y]= (a\cdot x)\cdot y-(-1)^{|x||y|}(a\cdot y)\cdot x \]
for $x,y\in \fg,$ $a\in M$ and $i,j\in\bbbz_2.$
We also  say {that}  $M$ together with a bilinear map $*\colon \fg\times M\longrightarrow M$ is a {\it left $\fg$-module} if
\[ \fg_i* M_j\sub M_{i+j} \quad \mbox{ and } \quad
[x,y]*a= x*(y* a)-(-1)^{|x||y|}y*(x*a)\]
 for  $x,y\in \fg,$ $a\in M$ and $i,j\in\bbbz_2.$
}
\end{deft}

\begin{rem}
{\rm
 We  note that if $(M,\cdot)$ is a right $\fg$-module, then  $M,$ together with the action
 \begin{equation}\label{inter}x* a:=-(-1)^{|x||a|}a\cdot x\;\;\;\;\;\hbox{ for $x\in \fg,a\in M$}
  \end{equation}is a left  $\fg$-module. Conversely, if $(M,*)$ is a left $\fg$-module, then  $M$ together with the action $ a\cdot x:=-(-1)^{|x||a|}x* a$ ($x\in \fg,a\in M$) is a right  $\fg$-module.}
  \end{rem}

{Although, if we have  a left $\fg$-module, we automatically have a right $\fg$-module and vice versa, our preference {is to}  use right actions
  as they  simplify working with degrees.}

 A linear map $\varphi$ from a $\fg$-module $M$ to a $\fg$-module $N$ is called a {\it $\fg$-module homomorphism} if $\varphi(mx)=\varphi(m)x$ for all $m\in M$ and $x\in \fg.$

 For a   $\fg$-module $M,$ a {\it derivation} of $\fg$ in $M$ is  a  linear map $d\colon \fg\longrightarrow M$  satisfying $$d[x,y]=d(x)y-(-1)^{|x||y|}d(y)x$$ for all $x,y\in \fg.$ We denote the set of all  derivations of $\fg$ in $M$ by  $\der(\fg,M).$ A~derivation $d\in \der(\fg,M)$ is called {\it inner} if there is $m\in M$ with $d(x)=mx$ for all $x\in\fg.$  The {\it first cohomology} $H^1(\fg,M)$  of $\fg$  with coefficients in $M$ is the quotient space {$\der(\fg, M)/{\rm Ider}(\fg,M),$  where}  ${\rm Ider}(\fg,M)$ is the set of inner derivations of $\fg$ in $M.$ A derivation of $\fg$ in $M$ is called {\it outer} if it is not inner.

\begin{deft}
{\rm Suppose that $\fk$ is a Lie superalgebra.

(i) For a  superspace $M,$ a bilinear  map $\alpha\colon \fk\times \fk\longrightarrow M$ is called  {\it supersymmetric} (resp. {\it skew-supersymmetric}) if for $x,y\in \fk,$
$\a(x, y)$ equals $(-1)^{|x||y|}\a(y, x)$ (resp. $-(-1)^{|x||y|}\a(y,x)$)
and it is called  {\it invariant} if
$$\a([x,y], z)=\a(x, [y,z]) \quad\quad       \hbox{for $x,y,z\in\fk.$}$$
The set of {all bilinear maps from $\fk\times\fk$ to $M$ is denoted by ${\rm Bil}_\bbbf(\fk,M)$} and the set of all supersymmetric invariant bilinear maps from $\fk\times \fk$ to $M$ is denoted by $\Sym(\fk, M)^\fk.$

 (ii) The subsuperspace
\[ \cent_{\bbbf}(\fk):=\lbrace \gamma \in \End(\fk)  |
(\forall a,b \in \fk)\ \gamma[a,b]=[\gamma(a),b]\rbrace\]
 of the superspace $\End(\fk)$ is called the {\it centroid} of $\fk.$
The Lie superalgebra   $\fk$ is called  {\it absolutely  simple} if $\frak{cen}_\bbbf(\fk)=\bbbf \hbox{id}.$
}
\end{deft}

Suppose that $\fk$ is a  finite dimensional Lie superalgebra and  $\kappa$ is an invariant nondegenerate supersymmetric bilinear form. Then the map
\begin{equation}\label{identify}
\varphi\colon \fk\longrightarrow \fk^*, \quad
\varphi(x) = \kappa(x,\cdot)
\end{equation}
is a linear bijection and so, for  $S\in\End(\fk),$
there is  a unique endomorphism $S^{*}$ of $\fk$ satisfying
\[ \kappa(Sx,y)=(-1)^{|x||y|}\kappa(S^{*}y,x) \quad \mbox{ for } \quad
x,y \in \fk.\]

\begin{lem}\label{basic}
For $T\in \End(\fk),$ define the bilinear map
{\begin{align*}
\kappa_T\colon & \fk\times \fk\longrightarrow \bbbf,\quad\quad
 (x, y)\mapsto \kappa(T(x),y).
\end{align*}}
This assignment has the following properties:
\begin{itemize}
\item[\rm(i)] $\kappa_T$ is supersymmetric (resp. skew-supersymmetric) if and only if $T^{*}=T$ (resp. $T^{*}=-T$).
\item[\rm(ii)] $\kappa_T$ is invariant if and only if $T\in\cent_\bbbf(\fk).$
\item[\rm(iii)] $\kappa_T$ satisfies $\kappa_T([x,y], z)=\kappa_T(x, [y,z])-(-1)^{|x||y|}\kappa_T(y,[x,z]),$ for all $x,y,z\in\fk,$ if and only if $T$ is a derivation.
\end{itemize}
\end{lem}

\begin{Pro}\label{cor-base}
Suppose that $\kappa$ is homogeneous and define
\[ {\theta\colon \End(\fk)\longrightarrow {\rm Bil}_\bbbf(\fk,\bbbf),
\quad\quad T\mapsto \kappa_{_T}.}\]
This assignment has the following properties:
\begin{itemize}
\item[\rm(a)] The restriction of $\theta$ to $$\der_-(\fk):=\{D\in \der(\fk)\mid D^{*}=-D\}$$ is a linear isomorphism from  $\der_-(\fk)$ onto the superspace $Z^2(\fk):=Z^2(\fk,\bbbf).$
\item[\rm(b)] The restriction of $\theta$ to the  superspace
$$\cent_+(\fk):=\{S\in \cent_\bbbf(\fk)\mid S^{*}=S\}$$  is a linear isomorphism from $\cent_+(\fk)$  onto the superspace $\Sym(\fk)^\fk:=\Sym(\fk,\bbbf)^\fk.$
\end{itemize}
\end{Pro}

\begin{deft}{\rm
The bilinear form $\kappa$ is called  {\it derivation invariant}  if  $\der(\fk)=\der_-(\fk).$}
\end{deft}
\begin{exa}
{\rm If $\fk$ is a sub-superalgebra of a Lie superalgebra $\fg$ such that $\kappa$ is a restriction of  an invariant supersymmetric bilinear form on $\fg$ and   each derivation of $\fk$ is of the form ${\rm ad}_x$ for some $x\in \fg,$ then $\kappa$ is derivation invariant.
}
\end{exa}
\begin{rem}\label{rem2}
{\rm
For $x,y\in\fk$ and $f\in \fk^*,$ one can define {$ (f\cdot x)(y):=f([x,y]).$}   Now $\fk^*$ together with this action is a $\fk$-module and $\varphi$  defined in (\ref{identify}) is a $\fk$-module isomorphism. Also for $\a\in Z^2(\fk),$ the linear map $\zeta_\a:\fk\longrightarrow \fk^*$ defined by $\zeta_\a(x)(y):=\a(x, y),$ for $x,y\in\fk,$  is  an element of $\der(\fk,\fk^*).$ Moreover, $\zeta_\a$ is an inner derivation if and only if $\a$ is a $2$-coboundary. Identifying $\fk$ and $\fk^*$ via $\varphi$ and using
Proposition~\ref{cor-base}, we  can embed {$H^2(\fk):=H^2(\fk,\bbbf)$  in $H^1(\fk,\fk)$ by considering $H^2(\fk)$} as  those outer derivations of $\fk$ belonging to $\der_-(\fk).$
}
\end{rem}

{{\bf Assumption:}
{\it From now on {to} the end of this section, we assume $\fk$ is a finite dimensional perfect Lie superalgebra equipped with a nondegenerate homogeneous invariant supersymmetric bilinear form $\kappa.$}}

\medskip
  {We set  $I$ to be the subsuperspace  of the exterior algebra\footnote{{For a vector superspace $V,$ by  $\Lam V=\op_{n=0}^\infty\Lam^nV$ and $S V=\op_{n=0}^\infty S^nV,$  we denote respectively the exterior superalgebra as well as the symmetric superalgebra  of the vector superspace $V$ and { ``$\;\wedge\;$'' and ``$\;\vee\;$'' denote  the multiplication maps on $\Lam~V$ and $S~V$ respectively.}}} $\Lam\fk$
     spanned by  $$[x,y]\wedge z-x\wedge[y,z]+(-1)^{|x||y|}y\wedge[x,z]\quad\quad \hbox{for  $x,y,z\in \fk.$}$$  Then}   the dual space  of the quotient space  $\Lam_{d}(\fk):=\Lam^2\fk/I$  is nothing but the superspace $Z^2(\fk)$ of  2-cocycles of $\fk$ with trivial coefficients.
Throughout this section, we fix  $\{D_1,\ldots, D_n\}\sub \der_-(\fk)$
such that $\{\kappa_{D_i}\mid 1\leq i\leq n\}$ is  a basis for $Z^2(\fk)$
(Proposition~\ref{cor-base}).
Since $\fk$ is finite dimensional, it follows that there is a  unique basis  $\{\lam_1,\ldots,\lam_n\}$ for  $\Lam_d(\fk)$ such that
\begin{equation}\label{alt}
\overline{x\wedge y}
=\sum_{i=1}^n \kappa_{D_i}(x,y)\lam_i
=\sum_{i=1}^n \kappa(D_ix,y)\lam_i\quad\quad\quad(x,y\in\fk),
\end{equation} where
$``\;\overline{\;\cdot\;}\;"$ stands for the equivalence classes in $\Lam_d(\fk).$
The degree 2-subspace   $S^2(\fk)$   has a natural $\fk$-module action and the dual space of the quotient space $S^2(\fk)/ S^2(\fk)\fk$ is isomorphic to $\Sym(\fk)^\fk.$  Now we fix a subset $\{S_1,\ldots,S_m\}$ of $\cent_+(\fk)$ such that $\{\kappa_{S_i}\mid 1\leq i\leq m\}$ is  a basis for $\Sym(\fk)^\fk.$ So there is a unique basis $\{\mu_1,\ldots,\mu_m\},$  which we fix throughout this section, for $S^2(\fk)/ S^2(\fk)\fk$ such that \begin{equation}\label{symm}
\overline{x\vee y}=\sum_{i=1}^m \kappa(S_ix,y)\mu_i\quad\quad\quad(x,y\in\fk),
\end{equation} in  which  by the abuse of notations, we again  use $``\;\overline{\;\cdot\;}\;"$   for the equivalence classes in $S^2(\fk)/ S^2(\fk)\fk.$
\medskip

Suppose $M$ is a superspace. If  $\a$ is a 2-cocycle of $\fk$ with coefficients in  $M,$ $\a$ induces the linear map
  \begin{align}\label{alt1}
{\widetilde \a}\colon &\Lam_d(\fk)\longrightarrow M, \quad\quad\overline{x\wedge y}\mapsto \a(x, y).
\end{align}
Next suppose  $\a\colon \fk\times \fk\longrightarrow M$ is a  supersymmetric invariant  bilinear map, then   $\a$ induces the linear map
  \begin{align}\label{hoch}
{\widehat \a}\colon &S^2(\fk)/ S^2(\fk)\fk\longrightarrow M,\quad\quad\overline{x\vee y}\mapsto \a(x, y).
\end{align}

\begin{Pro} \label{cocycle}
Let $M$ be a superspace. For $T\in\End(\fk)$ and $m\in M,$ define
\begin{align*}
   \nu_{_{T;m}}\colon &\fk\times \fk\longrightarrow M, \quad\quad(x, y)\mapsto \kappa_{_T}(x, y) m=\kappa(T(x),y)m.
\end{align*}
\begin{itemize}
\item[\rm(i)] If $m\in M$ and $S\in \cent_+(\fk),$ then $\nu_{_{S;m}}$ is  {a supersymmetric   invariant bilinear map. If } $D\in \der_-(\fk),$ then $\nu_{_{D;m}}$ is  a 2-cocycle of $\fk$ and  if $m\neq 0$, then
$\nu_{_{D;m}}$ is  a 2-coboundary if and only if $D$ is an inner derivation.
\item[\rm(ii)] If $\a\in Z^2(\fk,M),$ then $\a=\displaystyle{\sum_{i=1}^n\nu_{_{D_i;{\widetilde\a}(\lam_i)}}}$.
% KH: redundant because of (i).
%  {and} $\a$ is cohomologous to {{$\displaystyle{\sum_{i\in J}}\nu_{_{D_i;{\widetilde\a}%(\lam_i)}},$} where $J\sub\{1,\ldots,n\}$ is the set of those $i$  such that $D_i$ is an% outer derivation.
\item[\rm(iii)]
If $\a\colon \fk\times\fk\longrightarrow M$ is a supersymmetric invariant bilinear map, then \break $\a=\displaystyle{\sum_{i=1}^m\nu_{_{S_i;{\widehat\a}(\mu_i)}}}.$
\end{itemize}
\end{Pro}

\begin{prf}
(i) {By  Lemma \ref{basic},}  for $m\in M,$ $D\in \der_-(\fk)$ and $S\in \cent_+(\fk),$ $\nu_{_{D;m}}$ is a 2-cocycle and $\nu_{_{S;m}}$ is a supersymmetric invariant bilinear  map.  For the last statement, suppose  {that  $\nu_{_{D;m}}$ is a  2-coboundary. Then} there is a linear map $\ell\colon \fk\longrightarrow M$ such that for $x,y\in \fk,$ {$\kappa(Dx,y)m=\ell[x,y].$ This gives  a linear map $f\in\fk^*$  and a unique $t_f\in \fk$ such that for $x,y\in\fk,$  $$\kappa(Dx,y)=f[x,y]=\kappa(t_f,[x,y])=\kappa([t_f,x],y).$$ This in turn implies that $D={\rm ad}(t_f),$} in other words $D$ is an inner derivation. Conversely if $D$ is an inner derivation, it is immediate  that $\kappa_{_D}\fm m$ is  a coboundary as  {$\kappa$ is invariant.}

(ii) Considering (\ref{alt}) and (\ref{alt1}), for $x,y\in \fk,$ we have
{\small \begin{align*}
\a(x, y)={\widetilde\a}(\overline{x\wedge y})={\widetilde\a\Big(\sum_{i=1}^n \kappa(D_ix,y)\lam_i\Big)}=\sum_{i=1}^n \kappa(D_ix,y){\widetilde\a}(\lam_i)=\sum_{i=1}^n \nu_{_{D_i;{\widetilde\a}(\lam_i)}}(x, y).
\end{align*}}
This completes the proof.

(iii) Use the same argument as in part (ii).
\end{prf}

\subsection{$2$-cocycles of current superalgebras}
\label{subsec:4.2}

Throughout this subsection, $A$ denotes  a   unital supercommutative associative superalgebra.
For a superspace $M,$ we refer to a bilinear map  $F\colon A\times A\longrightarrow M$ satisfying
\begin{itemize}
\item[{(1)}] $F(a, b)=-(-1)^{|a||b|}F(b ,  a),$
\item[{(2)}] $F(ab ,  c)=F(a, bc)+(-1)^{|b||a|}F(b ,  ac) $
\end{itemize}
for $a,b,c\in A,$ a {\it Hochschild map}. For a Hochschild map $F\colon A\times A\longrightarrow M$ and $a\in A,$ using (2),  we have
{\begin{equation}\label{Hochschild}
F(1,a)=F(1,a1)=F(a,1)-F(a,1)=0.
\end{equation}}

\begin{deft}{\rm
  We set $\fg:=A\ot\fk$ and, for the sake of simplicity, for $a\in A$ and $x\in\fk,$ we denote $a\ot x$ by  $ax.$  We   recall  that $\fg$ together with  $$[ax,by]:=(-1)^{|x||b|}(ab)[x,y]$$ for $a,b\in A$ and $x,y\in \fk,$ is  a Lie superalgebra.}
\end{deft}

{\begin{deft}{\rm
Suppose  $\omega\colon \fg\times \fg\longrightarrow  M$ is a 2-cocycle of $\fg$ with coefficients in a superspace $M.$  For an element $a\in A$ and a homogeneous element $b\in A,$ define
\begin{align*}
\omega_{a,b}\colon &\fk\times\fk\longrightarrow M,\quad\quad (x, y)\mapsto (-1)^{|x||b|}\omega(a x,b y).
\end{align*}
We say $\omega$ is {\it $\fk$-cocyclic} if, for all homogeneous elements
$a,b\in A$ the form $\omega_{a,b}$ is a 2-cocycle.
We also  say $\omega$ is {\it $\fk$-invariant} if, for  all $a,b\in A,$ $\omega_{a,b}$ is an invariant bilinear map.
}
\end{deft}

\begin{exa}\label{exa-2-cocycle}
{\rm Suppose that $M$ is  a superspace. For a linear map $f\colon A\longrightarrow M$ and  $D\in \der_-(\fk),$
\[ \eta_{_{f,D}}\colon \fg\times  \fg\longrightarrow M, \quad
\eta_{f,D}(ax, by) := (-1)^{|b||x|}f(ab)\kappa(Dx,y),\;\;\;\;(a,b\in A,\;x,y\in\fk)\]
is a $\fk$-cocyclic $2$-cocycle of $\fg.$ Also for a Hochschild  map $F\colon A\times A\longrightarrow M$ and an element $S\in \cent_+(\fk),$
\[ \xi_{_{F,S}}\colon \fg\times \fg\longrightarrow M, \quad
\xi_{F,S}(ax, by) := (-1)^{|b||x|}F(a, b)\kappa(Sx,y),\;\;\;\;(a,b\in A,\;x,y\in\fk).\]
is a $\fk$-invariant $2$-cocycle of $\fg.$
}
\end{exa}}

\begin{lem}\label{basic1}
{\rm(i)} If  $\omega$ is  a $\fk$-cocyclic 2-cocycle, for  homogeneous elements $a,b\in A,$ we have
\[ \omega_{a,b}=(-1)^{|a||b|}\omega_{b,a}\andd
\omega_{ab,1}=\omega_{a,b} = \omega_{1,ab}. \]

{\rm (ii)} If  $\omega$ is  a $\fk$-invariant $2$-cocycle, then for homogeneous elements $a,b,c\in A,$ we have
{$$
\omega_{a,b}=-(-1)^{|a||b|} \omega_{b,a}\andd
\omega_{ab,c}=\omega_{a,bc}+(-1)^{|a||b|}\omega_{b,ac}.$$
}
\end{lem}

\begin{prf}
(i) Suppose that $a,b\in A$ and $x,y,z\in \fk$ are homogeneous elements. Then  we have
{\small \begin{align*}
\omega_{a,b}(x, y)=-(-1)^{|x||y|}{\omega}_{a,b}(y,x)= -(-1)^{|x||y|+|b||y|}\omega(ay,bx)=& (-1)^{|a||b|+|a||x|}\omega(bx,ay)\\
=& (-1)^{|a||b|}\omega_{b,a}(x,y).
\end{align*}}
We also have
{\small \begin{align*}
\omega_{ab,1}([y,z],x)&=\omega(ab[y,z],x)
=(-1)^{|y||b|}\omega([a y,b z], x)\\
&=(-1)^{|y||b|}\omega(a y,b [z, x])-(-1)^{|a||b|+|a||z|+|y||z|}\omega(b z,a [y,x])\\
&=\omega_{a,b}( y, [z, x])-(-1)^{|a||b|+|y||z|}\omega_{b,a}( z, [y,x])\\
&=\omega_{a,b}( y, [z, x])-(-1)^{|y||z|}\omega_{a,b}( z, [y,x])\\
&=\omega_{a,b}( [y, z], x).
\end{align*}}
This completes the proof as $\fk$ is perfect.
\smallskip

 (ii) Suppose that   $\omega$ is  a $\fk$-invariant  2-cocycle and $a,b\in A$ are {homogeneous. As $\omega_{a,b}$ is invariant, for $x,y,z\in\fk,$ we have
{\small\begin{align*}
\omega_{a,b}([x,y],z)&=\omega_{a,b}(x,[y,z])=-(-1)^{|y||z|}\omega_{a,b}(x,[z,y])\\
&=-(-1)^{|y||z|}\omega_{a,b}([x,z],y)=(-1)^{|y||z|+|x||z|}\omega_{a,b}(z,[x,y]).
\end{align*}}}
This shows that  $\omega_{a,b}$ is supersymmetric as $\fk$ is perfect. So for  $x,y,z\in\fk,$ we have
{\small\begin{align*}
\omega_{a,b}(x,y)=&(-1)^{|x||y|}\omega_{a,b}(y,x)=(-1)^{|b||y|+|x||y|}\omega(a y,b x)\\
=&-(-1)^{|a||b|+|a||x|}\omega(b x,a y)
=-(-1)^{|a||b|}\omega_{b,a}(x,y)
\end{align*}}
and

\begin{align*}
&\omega_{a,bc}(x,[y,z])+(-1)^{|a||b|}\omega_{b,ac}(x,[y,z])\\
=&(-1)^{|x||b|+|x||c|}\omega(ax,bc[y,z])+(-1)^{|a||b|+|x||a|+|x||c|}\omega(bx,ac[y,z])\\
=&(-1)^{|x||b|+|x||c|+|y||c|}\omega(ax,[by,cz])+
(-1)^{|a||b|+|x||a|+|x||c|+|c||y|}\omega(bx,[ay,cz])\\
=&(-1)^{|x||b|+|x||c|+|y||c|}\omega([ax,by],cz)+(-1)^{|x||c|+|y||c|+|a||b|+|a||y|+|y||x|}\omega(by,[ax,cz])\\
+&(-1)^{|a||b|+|x||a|+|x||c|+|c||y|}\omega(bx,[ay,cz])\\
=&(-1)^{|x||c|+|y||c|}\omega(ab[x,y],cz)+(-1)^{|y||c|+|a||b|+|a||y|+|y||x|}\omega(by,ac[x,z])\\
+&(-1)^{|a||b|+|x||a|+|x||c|}\omega(bx,ac[y,z])\\
=&\omega_{ab,c}([x,y],z)+(-1)^{|a||b|+|y||x|}\omega_{b,ac}(y,[x,z])+(-1)^{|a||b|}\omega_{b,ac}(x,[y,z])\\
=&\omega_{ab,c}([x,y],z)+(-1)^{|a||b|+|y||x|}\omega_{b,ac}(y,[x,z])-(-1)^{|a||b|+|x||y|}\omega_{b,ac}([y,x],z)\\
=&\omega_{ab,c}([x,y],z).
\end{align*}
This completes the proof as $\fk$ is perfect.
\end{prf}

\begin{Pro}\label{characterization1}
{ Suppose that $\fk$ is a finite dimensional perfect Lie superalgebra equipped with a nondegenerate invariant supersymmetric bilinear form $\kappa.$ If $\kappa$ is derivation  invariant}, then each $2$-cocycle of {$\fg=A\ot\fk$ is of the form $\omega_1+\omega_2,$} where $\omega_1$ is a $\fk$-invariant 2-cocycle and $\omega_2$ is a $\fk$-cocyclic $2$-cocycle.
\end{Pro}

\begin{prf}
Suppose that $\omega\colon \fg\times\fg\longrightarrow M$ is a $2$-cocycle. Then $\omega$ induces a linear map $f\colon \Lam^2\fg\longrightarrow M.$ The linear  maps
{\small
\begin{flalign*}
\sg_+\colon &\Lam^2 A\ot S^2\fk\longrightarrow \Lam^2\fg\\
& a\wedge b\ot x\vee y\mapsto \frac{1}{2}( (-1)^{|x||b|}(ax\wedge by)+(-1)^{|y||b|+|y||x|} (ay\wedge bx))\;\;\hbox{for $a,b\in A,$ $x,y\in \fk$}
\end{flalign*}}
and
{\small\begin{flalign*}
\sg_-\colon &S^2 A\ot \Lam^2\fk\longrightarrow \Lam^2\fg\\
& a\vee b\ot x\wedge y\mapsto \frac{1}{2}( (-1)^{|x||b|}(ax\wedge by)-(-1)^{|y||b|+|y||x|} (ay\wedge bx))\;\; \hbox{for $a,b\in A,$ $x,y\in \fk$}
\end{flalign*}}
are embeddings and $$\Lam^2\fg\simeq (\Lam^2 A\ot S^2\fk)\op  (S^2 A\ot \Lam^2\fk).$$
On the other hand, the linear  map
\[ \begin{array}{rl}
e_1\colon A\longrightarrow S^2A,\quad\quad a\mapsto a\vee 1,
\end{array}\]
is also an embedding. The linear  map
\begin{align*}
\mu\colon S^2 A\longrightarrow A,\quad\quad
a\vee b\mapsto ab,
\end{align*}
satisfies $S^2A= {\rm im}(e_1)\op {\rm ker}(\mu).$ So
\[ \Lam^2\fg\simeq (\Lam^2 A\ot S^2\fk)\op
( A\ot \Lam^2\fk)\op ({\rm ker}(\mu)\ot\Lam^2 \fk).\]
 We identify $f$  as $f_1+f_2+f_3$ where
\[ f_1\colon \Lam^2 A\ot S^2\fk\longrightarrow M,\;\; f_2\colon  A\ot
\Lam^2\fk\longrightarrow M,\;\; f_3\colon {\rm ker}(\mu)\ot\Lam^2 \fk\longrightarrow M.\]
In view of Remark \ref{rem2} and the fact that $\fk$ is perfect and
$\kappa$ is derivation invariant,
the super versions of {Corollary 3.3, Theorem 3.7 and (13)} of \cite{NW},
imply  that $f_3=0$ and that $f_1$ and $f_2$ respectively induce maps
$\widetilde{f_1}\colon A\times A\longrightarrow \Sym(\fk,M)^\fk$ and
$\widetilde{f_2}\colon A\longrightarrow Z^2(\fk,M)$ with \[\tilde{f_1}(a,b)(x, y)=f_1((a\wedge b)\ot (x\vee y))\andd \tilde{f_2}(a)(x, y)=f_2(a\ot (x\wedge y)).\]
We next define
$\omega_j\colon \fg\times \fg\longrightarrow M$, $j=1,2$, by
\[ \omega_1(ax, by) :=  (-1)^{|x||b|}f_1(a\wedge b\ot x\vee y)
\quad\quad \hbox{for $a,b\in A$ and $x,y\in \fk$} \]
and
\[ \omega_2(ax, by) :=  (-1)^{|x||b|}f_2(a b\vee 1\ot x\wedge y)\quad\quad
\hbox{for $a,b\in A$ and $x,y\in \fk.$} \]
Then  $\omega_1$ and $\omega_2$ are $2$-cocycles (see \cite[Thm.~3.7]{NW}). { Also  $f=f_1+f_2$ implies $\omega=\omega_1+\omega_2.$} As ${\rm im}(\widetilde f_2)\sub \Sym(\fk,M)^\fk,$ $\omega_1$ is a $\fk$-invariant $2$-cocycle and as ${\rm im}(\widetilde{f_2})\sub Z^2(\fk, M),$ $\omega_2$ is a $\fk$-cocyclic $2$-cocycle.
\end{prf}

\begin{Pro}\label{characterization2}
{For $D_i,\lam_i$ and $S_j, \mu_j$} from \eqref{alt} and \eqref{symm},
the following assertions hold for any
$\omega \in Z^2(\g,M)$:
\begin{itemize}
\item[\rm(i)] If $\omega$ is $\fk$-cocyclic,
then  there are linear maps $f_1,\ldots,f_n\in{\rm Hom}_\bbbf( A,M)$ such that
$\omega = \sum_{i=1}^n \eta_{f_i, D_i}$.
\item[\rm(ii)] If $\omega$ is  $\fk$-invariant,
then there are  Hochschild  maps $F_i\colon A\times A\longrightarrow M,
i=1,\ldots,m$, such that $\omega = \sum_{i=1}^m \xi_{F_i, S_i}$.
\end{itemize}
\end{Pro}

\begin{prf}
  (i) For  $a,b\in A,$ $\omega_{a,b}$ is a 2-cocycle. Consider ${\widetilde\omega}_{a,b}$ as in (\ref{alt1}) and for $i\in\{1,\ldots,n\},$  take
\begin{align*}
f_i\colon &A\longrightarrow M, \qquad
f_i(a) := {\widetilde\omega}_{1,a}(\lam_i).
\end{align*}
Then $f_i$ is  a linear map and by  Proposition \ref{cocycle}(iii)
and Lemma \ref{basic1}(i), we have
\begin{align*}
\omega_{a,b}=\sum_{i=1}^n \nu_{_{D_i;{\widetilde\omega}_{a,b}(\lam_i)}}=\sum_{i=1}^n \nu_{_{D_i;{\widetilde\omega}_{1,ab}(\lam_i)}}=\sum_{i=1}^n \kappa_{D_i} f_i(ab).
\end{align*}

(ii) For  $a,b\in A,$  consider ${\widehat\omega}_{a,b}$ as in (\ref{hoch}) and for $i\in\{1,\ldots,m\},$  take
{\begin{align*}
F_i\colon &A\times A\longrightarrow M,\quad\quad(a, b)\mapsto {\widehat\omega}_{a,b}(\mu_i).
\end{align*}}
Then by Lemma \ref{basic1}, $F_i$ is a Hochschild map and  by {Proposition} \ref{cocycle}, we have
\begin{align*}
\omega_{a,b}=\sum_{i=1}^m \nu_{_{S_i;{\widehat\omega}_{a,b}(\mu_i)}}=\sum_{i=1}^m  F_i(a, b)\kappa_{S_i}.
\end{align*}This completes the proof.
\end{prf}

 Recalling Example \ref{exa-2-cocycle} and using Propositions \ref{characterization1} and \ref{characterization2}, we arrive at the following structure theorem
for $2$-cocycles:

\begin{Thm}\label{cor1}
 Suppose that $\kappa$ is derivation invariant and $D_1,\ldots, D_t \in
\der(\fk)$ are such that the cohomology classes $[\kappa_{D_j}]$
form a basis of $H^2(\fk)$ and
$S_1, \ldots, S_m \in \cent_+(\fk)$ are such that the
$\kappa_{S_i}$ form a basis of $\Sym(\fk)^\fk.$
Then  each 2-cocycle of $\fg=A\ot \fk$  with values in a superspace $M,$  is cohomologous to a sum
\[ \tilde \omega = \sum_{i=1}^t \eta_{_{f_i,D_i}}+\sum_{j=1}^m\xi_{_{F_j,S_j}}, \]
i.e.,
\[ \tilde\omega(ax,by)
= (-1)^{|b||x|} \sum_{i=1}^t f_i(ab)\kappa(D_i x,y)
+ (-1)^{|b||x|} \sum_{j=1}^m F(a, b)\kappa(S_j x,y),\]
where the $f_i \: A \to M$ are  linear maps
and the $F_j\colon A\times A\longrightarrow M$ are Hochschild  maps.
\end{Thm}

\section{Unitary representations of current superalgebras}
\label{sec:5}

A finite dimensional real Lie algebra $\fg$ is called {\it compact}
if it is the Lie algebra of a compact Lie group~$G.$
Then the  subgroup of ${\rm Aut}(\fg)$ generated by ${\rm e}^{{\rm ad}(\fg)}$
is compact {\cite[Pro. 12.1.4]{HN}}.
Further, compactness of $\fg$ is equivalent to the existence of a
faithful finite dimensional unitary representation
\cite[Thm.~12.3.9, Lem.~12.1.2]{HN}.
This is different for Lie superalgebras.
We call a finite dimensional real  Lie superalgebra
$\fg=\fg_0\op\fg_1$ {\it compact}
if the subgroup of ${\rm Aut}(\fg)$ generated by {${\rm e}^{{\rm ad}(\fg_0)}$}
has compact closure. As we have already seen in the introduction,
this does not imply the
existence of a faithful finite dimensional unitary representation.
%(see Theorem~\ref{S-K}, Remark~\ref{rem4}(iii) and Lemma \ref{necessary}).
For a classification of compact Lie superalgebras
with faithful unitary representations we refer to~\cite{AN}.
Our aim in this  section is {to  investigate} the existence of
(projective)
unitary representations for current superalgebras of the form
$ A\ot \fk$, where $\fk$ is a  simple compact Lie superalgebra
and $A$ is graded commutative.

\subsection{Compact Lie algebras}
\label{subsec:4.3}

We start with the case {where} $\fk$ is a simple compact Lie algebra.
Fix a simple compact  Lie algebra $\fk$  with the Killing form $\kappa$ and
a compact Lie group~$K$ with Lie algebra $\fk$.

\begin{rem}\label{rem3}
{\rm  Suppose that $A$ is a unital supercommutative associative superalgebra.
Since  $\fk$  is a  compact simple Lie algebra, {$H^2(\fk)=\{0\}$ and $\cent_\bbbr(\fk)=\bbbr{\rm id}.$ By Theorem~\ref{cor1}, a} $2$-cocycle $\omega$  of  $A\ot\fk$
is equivalent to one of the form
\begin{equation}\label{ref}
\omega(ax, by):= \omega_F(ax, by)
:= F(a,b)\kappa(x,y)\quad\quad \hbox{for $a,b\in A,\; x,y\in \fk$},
\end{equation}
where $F$ is a  Hochschild  map. Here we use that $|x| = 0$ for every $x \in
\fk = \fk_0$. }
\end{rem}

In the following proposition we use Theorem~\ref{main1} from the appendix.

\begin{Pro}\label{unitary-graded}
{Let $M$ be a superspace and let $A=A_0\op A_1$ be
a finite dimensional unital supercommutative  associative  superalgebra equipped with  a
consistent\footnote{This means that the $\bbbz_2$-grading induced from the
$\bbbz$-grading on $A$ coincides with the original $\bbbz_2$-grading on $A.$}
$\bbbz$-grading $A=A^0\op A^1\op A^2$ satisfying
$A^0=\bbbr$ and $A^1A^1=A^2.$
Consider the central extension $\hat\g$
of the Lie superalgebra $\fg:=A\ot \fk$
defined by a 2-cocycle of the form
$\omega = \omega_F$, where $F:A\times A\longrightarrow M$ is an
even Hochschild map.
Then $\fn:=((A^1\op A^2)\ot \fk)\op M$
is a Clifford--Lie superalgebra and
$\widehat\fg\simeq\fn\rtimes  \fk.$
If $K$ is a simply connected Lie group with Lie algebra $\fk$,
then the equivalence classes of irreducible unitary representations of the
Lie supergroup $\cN \rtimes K$
are determined by the $K$-orbits $\cO_\lambda$ in $\sc(\fn)^\star$
and, for any fixed $\lambda$, they are in one-to-one correspondence
with odd irreducible representations of $\hat K_\lambda^\circ;$ see (\ref{K-circle}).}
\end{Pro}

\begin{prf}
We identify $A^0\ot \fk$ with $\fk,$ so  that  $\fn$ is an ideal of $\widehat\fg$ and $\widehat\fg\simeq\fn\rtimes \fk.$  {As $A^2=A^1A^1$ and $A^3=\{0\},$ we get  $F(A^2,A^2)=\{0\}.$} We next note that for  $X\in A^2\ot\fk$ and   $Y\in (A^1\ot \fk)\op(A^2\ot \fk),$ we have $\omega(X,Y)=0$ as $F$ is even and  {$F(A^2,A^2)=\{0\}$ and so}
$$[X,Y]_\omega=[X,Y]+\omega(X,Y)=0.$$
{This implies  that {$\fn_0=(A^2\ot \fk)\op M_0\sub Z(\fn),$} i.e., $\fn$
is an ideal of $\widehat\fg$ which is a
Clifford--Lie superalgebra. Now the assertion follows from
Theorem~\ref{main1}.}
\end{prf}

\begin{Pro}\label{3-elm}
Assume $A=A_0\op A_1$ is a unital supercommutative  associative superalgebra and  $\omega$ is a $2$-cocycle for $\fg:=A\ot \fk$ with corresponding central extension  $(\widehat\fg,[\cdot,\cdot]_\omega).$ Suppose {$x\in \fk$  and  $a_1,\ldots,a_n\in A,$ $n\geq 2,$} are homogeneous elements  {with $a_1^2=\cdots=a_n^2=0.$ Then $X:= a_1\cdots a_n\ot x$ satisfies   $[X,X]_\omega=0.$} In particular, if $n\geq 3$ is odd and $a_1,\ldots,a_n$ are odd elements, $X\in {\rm urad}(\widehat\fg)$.
\end{Pro}

\begin{prf}
Recall Remark \ref{rem3} and suppose  $F$ is a Hochschild map  on $A\times A$ such that  for $a,b\in A$ and $x,y\in\fk,$ $\omega(ax,by)=F(a, b)\kappa(x,y).$
 { For $a:=a_1\cdots a_t,$   we have}
\begin{align*}
F(a, a)=F(a_1, a_2\cdots a_ta_1\cdots a_t)+(-1)^{|a_1|(|a_2|+\cdots+|a_t|)}F(a_2\cdots a_t, a_1a_1\cdots a_t)=0
\end{align*}
{For  $x\in \fk,$  this leads to} $$[ a_1\cdots a_t\ot x, a_1\cdots a_t\ot x]_\omega=\kappa(x,x)F(a_1\cdots a_t, a_1\cdots a_t)=\kappa(x,x)F(a, a)=0.$$ {If $t\geq 3$ is odd and $a_1,\ldots,a_t$ are odd elements, $X:=a_1\cdots a_t\ot x$ is an odd element of $\widehat\fg$ so that the remaining assertion follows from  Lemma \ref{necessary}.}
\end{prf}

\begin{exa}{\rm
A typical (and universal)
example of a unital supercommutative associative superalgebra is the real  unital supercommutative associative superalgebra $\Lam_s(\bbbr)$  generated by odd elements $\ep_1,\ldots,\ep_s$ subject to the relations $$\ep_i\ep_j+\ep_j\ep_i=0\quad\quad\quad(1\leq i,j\leq s).$$ {It is called the (real) {\it Gra\ss{}mann superalgebra} in $s$ generators}
$\ep_1,\ldots,\ep_s.$
The Gra\ss{}mann superalgebra $\Lam_s(\bbbr)$ has a natural consistent  $\bbbz$-grading $$\Lam_s(\bbbr)={\bigoplus_{m\in\bbbz}}\Lam_s^m(\bbbr)$$ with $$\Lam_s^0(\bbbr)=\bbbr,\;\Lam_s^m(\bbbr)=\hbox{span}_\bbbr\{\ep_{i_1}\cdots\ep_{i_m}\mid 1\leq i_1<\cdots<i_m\leq s\}\andd\Lam_s^n(\bbbr)=\{0\}$$
for $1\leq m\leq s$ and $n\in \bbbz\setminus\{0,\ldots,m\}.$ {We set
$$\Lam_s^{\rm odd}(\bbbr):=\bigoplus_{m=0}^\infty\Lam_s^{2m+1}(\bbbr),\;\;\Lam_s^{\rm even}(\bbbr):=\bigoplus_{m=1}^\infty\Lam_s^{2m}(\bbbr)\andd\Lam_s^{+}(\bbbr):=\bigoplus_{m=1}^\infty\Lam_s^{m}(\bbbr).$$
}
}\end{exa}

The following theorem reduces the problem of classifying projective irreducible unitary
representation of $\g = \Lambda_s(\R) \otimes \fk$ to the case of
semidirect products of $\fk$ with Clifford--Lie superalgebras
discussed in detail in Appendix~\ref{sec:3}.

\begin{Thm}\label{urad}
Consider a {current superalgebra} $\g = \Lambda_s(\R) \otimes \fk$ and
a central extension $\g = \g \oplus M$ by a superspace $M$, defined by
a {2-cocycle} $\omega = \omega_F$, where
$F\colon \Lam_s(\bbbr)\times \Lam_s(\bbbr)\longrightarrow M$ is a Hochschild map.
We set
\[ R:=\hbox{\rm span}_\bbbr\{F(a, b)\mid a\in \Lam_s(\bbbr),b\in \op_{m=3}^s\Lam_s^m(\bbbr)\}
\subeq M.\]
\begin{itemize}
\item[\rm(i)] The ideal
\[ I:=(\op_{m=3}^s\Lam_s^m(\bbbr)\ot \fk) \op R \]
 of $\widehat\fg$ lies in the unitary radical   ${\rm urad}(\widehat\fg)$  of $\widehat\fg;$ in particular,  for each   unitary representation $\pi$ of $\widehat\fg,$ we have  $I\sub{\rm ker}(\pi).$
\item[\rm(ii)] Suppose that $\omega$ is even, i.e., that $F$
is even.  Set
\[  \widehat\fn:=(\Lam^+(\bbbr)\ot \fk)\op M
\quad \mbox{ and } \quad
\fn:=\widehat\fn/I.\]
 Then $\widehat\fg\simeq \widehat\fn\rtimes \fk$ and $\widehat\fg/I\simeq \fn\rtimes \fk.$  Moreover, $\fn$ is a Clifford--Lie superalgebra and  if
$K$ is a simply connected Lie group with Lie algebra $\fk$,
then the equivalence classes of irreducible unitary representations of the
Lie supergroup $\cN \rtimes K$
are determined by the $K$-orbits $\cO_\lambda$ in $\sc(\fn)^\star$
and, for any fixed $\lambda$, they are in {one-to-one} correspondence
with odd irreducible representations of $\hat K_\lambda^\circ$ as defined in \eqref{K-circle}).
\end{itemize}
\end{Thm}

\begin{prf}
(i) { Assume  $r,t$ are positive integers  and  $j_1,\ldots,j_{2t+2},i_1,\ldots,i_r\in\{1,\ldots,s\}.$ We know from Proposition~\ref{3-elm} that
\begin{equation*}\label{in-ker0}
\ep_{j_1}\ldots\ep_{j_{2t+1}}\ot \fk\sub {{\rm urad}(\widehat\fg)}.
\end{equation*}
 For    $x\in \fk,$ we  have
\begin{fleqn}
\begin{align*}
F(\ep_{j_1}\cdots\ep_{j_{2t+1}}, \ep_{i_1}\cdots\ep_{i_r})\kappa(x,x)
&= [\ep_{j_1}\cdots\ep_{j_{2t+1}}\ot x,\ep_{i_1}\cdots\ep_{i_r}\ot x]_\omega\\
&\in [{\rm urad}(\widehat\fg),\widehat\fg]_\omega\sub{\rm urad}(\widehat\fg).
\end{align*}
\end{fleqn}
Now choosing $x\in\fk$ with $\kappa(x,x)\neq0,$ we get that
\begin{equation}\label{in-ker1}
F(\ep_{j_1}\cdots\ep_{j_{2t+1}}, \ep_{i_1}\cdots\ep_{i_r})\in {\rm urad}(\widehat\fg)
\end{equation}
This implies {that,} for $x,y\in\fk,$
\begin{align*}
&\ep_{j_1}\cdots\ep_{j_{2t+1}}\ep_{i_1}\cdots\ep_{i_r}\ot [x,y]\\
=&[\ep_{j_1}\cdots\ep_{j_{2t+1}}\ot x,\ep_{i_1}\cdots\ep_{i_r}\ot y]_\omega
-F(\ep_{j_1}\cdots\ep_{j_{2t+1}}, \ep_{i_1}\cdots\ep_{i_r})\kappa(x,y)\\
\in& [{\rm urad}(\widehat\fg),\widehat\fg]_\omega+{\rm urad}(\widehat\fg)\sub{\rm urad}(\widehat\fg).
\end{align*}
{As} $\fk$ is perfect, this implies that
\begin{equation}\label{in-ker2}
\op_{m=3}^s\Lam_s^m(\bbbr)\ot \fk\sub {\rm urad}(\widehat\fg).
\end{equation} This in turn {shows that,}  for $x\in \fk,$
\begin{fleqn}
\begin{align*}
F(\ep_{j_1}\cdots\ep_{j_{2t+2}}, \ep_{i_1}\cdots\ep_{i_r})\kappa(x,x)=[\ep_{j_1}\cdots\ep_{j_{2t+2}}\ot x,\ep_{i_1}\cdots\ep_{i_r}\ot x]_\omega\in {\rm urad}(\widehat\fg).
\end{align*}
\end{fleqn}
{Choosing} $x\in\fk$ with $\kappa(x,x)\neq 0,$ one gets that   $F(\ep_{j_1}\cdots\ep_{j_{2t+2}}, \ep_{i_1}\cdots\ep_{i_r})$ lies in ${\rm urad}(\widehat\fg).$  This together with  (\ref{in-ker1}) and the fact that $F(1, \Lam_s(\bbbr))=\{0\}$ gives that $R\sub {\rm urad}(\widehat\fg).$ So by (\ref{in-ker2}), we get that $$I=(\op_{m=3}^s\Lam_s^m(\bbbr)\ot \fk) \op R\sub {\rm urad}(\widehat\fg).$$ This completes the proof as ${\rm urad}(\widehat\fg)$ lies in the kernel of each unitary representation of $\widehat\fg$ by Lemma \ref{necessary}.}

(ii) Set $$A:=\Lam_s(\bbbr)/(\op_{m=3}^s\Lam_s^m(\bbbr))\andd T:=M/R$$ and consider the {induced Hochschild map  $\bar F:A\times  A\longrightarrow  T.$   Then $$\widehat\fg/I\simeq (A\ot \fk)\op T$$ is a central extension of $A\ot \fk$ corresponding to  the $2$-cocycle $\xi_{\bar F,{\rm id}}.$  Therefore, irreducible unitary representations of  $\widehat\fg$
are exactly those of $\widehat\fg/I\simeq(A\ot \fk)\op T\simeq\fn\rtimes\fk $ and so  the assertion  follows from Proposition \ref{unitary-graded}.}
\end{prf}

\begin{Pro}
The  Lie superalgebra $\fg:=\Lam_s(\bbbr)\ot_\bbbr \fk$
has a central extension with a faithful unitary representation if and only if
$s \leq 2$.
\end{Pro}

\begin{prf} Suppose first that $s\geq 3.$ For all  $x\in \fk$ and
distinct indices $i_1,\ldots, i_{2m+1}$ ($m\geq 1$),
Proposition \ref{3-elm} implies that $\ep_{i_1}\ldots\ep_{i_{2m+1}}\ot x$
 has square zero in each central extension $\hat\g$ of $\fg$.
In particular, $\hat\g$ has no faithful unitary representation.

Now we assume that $s\in\{1,2\}$. We consider the even  bilinear map
\begin{align*}
F\colon &\Lam_s(\bbbr)\times \Lam_s(\bbbr)\longrightarrow \bbbr\\
&(\ep_i,\ep_j)\mapsto \d_{i,j}\;\;\;\;\;( 1\leq i,j\leq s)\\
&(b,a),(a,b)\mapsto 0\;\;\;\;\; (a\in\Lam_s(\bbbr),\; b\in \Lam_s(\bbbr)_0).
\end{align*}
We shall show  that $F$ is a Hochschild map: Suppose that $a,b,c$ are homogeneous elements of $\Lam_s(\bbbr)$ with respect to the $\bbbz$-grading on $\Lam_s(\bbbr).$ If at least one of $a,b$ is even, then we have $F(a,b)=F(b,a)=0.$ If both $a$ and $b$ are odd, then we may assume $a=\ep_i$ and $b=\ep_j$ for some $1\leq i,j\leq s$ which in turn implies that $F(a,b)=F(b,a).$ Therefore in both cases
\begin{equation}\label{prop1}
F(a,b) =-(-1)^{|a||b|}F(b,a).
\end{equation}
Moreover, we know that  $F$ is even,   $F(\Lam_s(\bbbr)_0,\Lam_s(\bbbr))=F(\Lam_s(\bbbr),\Lam_s(\bbbr)_0)=\{0\}$ and $\Lam_s^3(\bbbr)=\{0\},$
so  we get that $F(ab,c)=F(a,bc)=F(b,ac)=0$. In fact,
the only critical case for $F(ab,c)$ is that, say $a \in \Lambda_s^0(\R)$ and
$b,c \in \Lambda^1_s(\R)$. Then
$F(a,bc)=0$ and by (\ref{prop1}), we have $F(ab,c)=F(b,ac).$ Therefore
we always have
\[ F(ab,c)=F(a,bc)+(-1)^{|a||b|}F(b,ac) \]
and thus $F$ is a Hochschild map. So \begin{align*}
\omega\colon \fg\times\fg\longrightarrow \bbbr,\quad\quad (ax,by)\mapsto F(a, b)\kappa(x,y)
\end{align*}
is a $2$-cocycle whose restriction to $\fg_1\times\fg_1$ is definite as  the Killing form $\kappa$ of   $\fk$ is negative definite ($\fk$ is compact and simple).
Finally Lemma~\ref{point-crit} implies that
the cone $\sc(\hat\g)$ of the central extension $\hat\g$ is pointed.
As $\hat\g \cong \fn \rtimes \fk$ is a semidirect product of the compact
Lie algebra $\fk$ with a Clifford--Lie superalgebra,
Theorem~\ref{main1} now implies that the unitary representations of
$\hat\g$ separate the points. As we may form arbitrary direct sums,
a faithful unitary representation exists.
\end{prf}

\subsection{Compact  simple Lie superalgebras}
\label{subsec:5.2}
By Theorem 2.3 of \cite{AN}, one knows the classification of finite dimensional compact Lie superalgebras; to state that classification, we first need to recall
the structure of the Lie superalgebras  involved.

Denote by  $M_{m\times n}(\bbbc)$  the set of all $m\times n$-matrices with complex entries and for  a complex matrix $A,$ {denote   by  $A^t$ and $A^*$   the transposition and the conjugation of the  matrix~$A$ respectively.}  For a square matrix (resp.  block matrix)~$A,$ by ${\rm tr}(A)$ (resp. ${\rm str}(A)$), we mean the trace (resp. supertrace) of~$A.$
  We set $$\fu(n;\bbbc):=\{A\in M_{n\times n}(\bbbc)\mid A^*=-A\}\hbox{~~and~~} \frak{su}(n;\bbbc):=\{A\in \fu(n;\bbbc)\mid {\rm str}(A)=0\}.$$
Also for  two positive integers~$p,q,$ we {denote by $\frak{gl}(p,q)$  the}  Lie superalgebra of all block matrices of dimension $(p,q)$ with entries in $\bbbc$ and {denote by  $\frak{sl}(p,q)$  the sub-superalgebra}  of $\frak{gl}(p,q)$ containing all elements with zero supertrace.   For  $X=\left(\begin{array}{cc}A&B\\C&D\end{array}\right)\in\frak{gl}(p,q),$ we put the {\it superconjugation} of $X$ to be $$X^{\#}:=\left(\begin{array}{cc}A^*&-iC^*\\-iB^*&D^*\end{array}\right)\in\frak{gl}(p,q)$$  and set
{\small \begin{align*}
\fu(p|q;\bbbc):=\{X\in \frak{gl}(p,q)\mid X^\#=-X\}=\left\{\left(\begin{array}{cc}A&B\\iB^*&D\end{array}\right)\mid A\in\fu(p;\bbbc), D\in \fu(q;\bbbc)\right\}
\end{align*}}
which  is a compact real form {of}  $\frak{gl}(p,q).$

\medskip

\noindent {\bf I}. $\pmb{\frak{su}(p|q;\bbbc)}$ {\bf and} $\pmb{\frak{psu}(p|p;\bbbc)}\colon $ Suppose $p,q$ are two positive integers with $p\geq q$ and let ${\bf 1}_p$ (resp. ${\bf 1}_q$) be the identity matrix of order $p$ (resp. $q$). We  set  $$\fsu(p|q;\bbbc):=\left\{X\in \fu(p|q;\bbbc)\mid \;{\rm str}(X)=0\right\} $$  and
$$\fpsu(p;\bbbc):=\fsu(p|p;\bbbc)/\bbbr i{\bf 1}. $$
$\frak{su}(p|q;\bbbc)_0\simeq \frak{su}(p;\bbbc)\op\frak{su}(q;\bbbc)\op i\bbbr\i$ with $\i:=\frac{1}{p}{\bf 1}_p+\frac{1}{q}{\bf 1}_q$ and $\frak{su}(p|q;\bbbc)_1\simeq \bbbc^p\ot \overline{\bbbc^q}.$ {Since $i\bbbr{\bf 1}$ has a trivial intersection with $\frak{su}(p|p;\bbbc)_1,$ to simplify our notations, we take $\frak{psu}(p|p;\bbbc)_1=\frak{su}(p|p;\bbbc)_1.$}
 The supertrace form $\kappa,$ that is the  bilinear form mapping $(A,B)$ to ${\rm str}(AB),$ is an even supersymmetric bilinear form on $\fsu(p|q;\bbbc)$ whose restriction to $\frak{su}(p;\bbbc)$ is negative definite while its restriction to $\frak{su}(q;\bbbc)$ is positive definite. Moreover, if $p>q,$ the restriction of $\kappa$ to  $ i\bbbr\i$ is positive definite while the  radical of $\kappa$  is  $\bbbr i{\bf 1}$ if $p=q;$ in the latter case, $\kappa$ induces a nondegenerate even supersymmetric bilinear form on $\frak{psu}(p|p;\bbbc)$ denoted again by $\kappa.$ These all together imply that
\begin{equation}\label{equalzero1}
\parbox{4in}{{in both cases $\frak{su}(p|q;\bbbc)$ ($p>q$) and $\frak{psu}(p|p;\bbbc),$   {there exist}   nonzero even elements   $x,y$ with
$0\neq \kappa(x,x)=-\kappa(y,y)$ and  $\kappa(x,y)=0.$}}
  \end{equation}

We now  suppose {that}  $p>q\geq 1$ and note that $\frak{su}(p|q;\bbbc)_1$ is an irreducible $\frak{su}(p|q;\bbbc)_0$-module  on which $ i\i$ acts as $(\frac{1}{p}-\frac{1}{q}){\rm id}.$   Also $\frak{psu}(p|p;\bbbc)_0$ {acts} irreducibly  on $\frak{psu}(p|p;\bbbc)_1\simeq \bbbc^p\ot \overline{\bbbc^p}.$

{The Lie superalgebra    $\fsu(p|q;\bbbc)$ is a compact real form of $\frak{sl}(p,q)$ and} $\fpsu(p;\bbbc)$ is a compact real form of $\frak{sl}(p,p)/\bbbc{\bf 1}.$
On the other hand, the Cartan--Killing form of $\fsu(p|q;\bbbc)$ is nondegenerate which  in turn implies that all derivations of $\frak{su}(p|q;\bbbc)$ are inner and so  $$H^2(\fsu(p|q;\bbbc))=\{0\}$$ while  the restriction of $D:={\rm ad}\left(\begin{array}{cc}
i{\bf 1}_p&0\\0&0
\end{array}\right)$ to $\frak{su}(p|p;\bbbc)$   induces  an even outer derivation of  $\frak{psu}(p|p;\bbbc)$ vanishing on the even part; this  in fact forms  a  basis for $H^2(\frak{psu}(p|p;\bbbc))\simeq\bbbr$ (see Remark \ref{rem2}).

For our future use, we note that
for
\begin{equation}\label{special}B:=\diag(1,-1,0\ldots,0),\;\;\;x_*:=\left(\begin{array}{cc}
0& B\\
iB& 0
\end{array}
\right), y_*:=\left(\begin{array}{cc}
0& {\bf 1}_n\\
i{\bf 1}_n& 0
\end{array}
\right)\in \frak{psu}(p|p;\bbbc)_1,
\end{equation} we have
\begin{equation}\label{f-use}
\kappa(x_*,y_*)=\kappa(Dx_*,y_*)=0\andd [x_*,y_*]=u+v,
\end{equation}
 for some nonzero elements $u$ and $v$ of irreducible components of $\frak{psu}(p|p;\bbbc)_0.$
Also for the elementary matrix $e_{r,s}$ with $1$ in $(r,s)$-entry and $0$ elsewhere and
\begin{equation}\label{z*}
z_*:=\left(\begin{array}{cc}
0& e_{r,s}\\
ie_{s,r}& 0
\end{array}
\right)\in \frak{su}(p|q;\bbbc)_1,
\end{equation} we have
\begin{equation}\label{uv}
\kappa(z_*,z_*)=0\andd [z_*,z_*]=u+v+z
\end{equation} where $u,v$ are nonzero elements of irreducible components of $\frak{su}(p|q;\bbbc)_0$ and $z$ is a central element of $\frak{su}(p|q;\bbbc)_0.$

 \noindent {\bf II}.  $\pmb{\frak{c}(n)}$ $\pmb{(n\geq 2)}:$ Suppose $n\geq 2$ and let $\frak{c}(n)$  be the
compact real form of the orthosymplectic Lie superalgebra  $\frak{osp}(2|2n-2)$ with {$\frak{c}(n)_0\simeq\bbbr\op\frak{sp}(n-1),$ where $\frak{sp}(n-1)$  is the compact real form of  the symplectic Lie algebra of rank $n-1$} and $\frak{c}(n)_1$  {is an irreducible $\frak{c}(n)_0$-module} isomorphic to $\mathbb{H}^{n-1}$ {where $\mathbb{H}$} is the  quaternion algebra; more precisely, $\frak{c}(n)$ is the set of all matrices
{\small $$\left(
\begin{array}{cccc}
\a &0& M&N\\
0&-\a&i\overline N&-i\overline M\\
-i\overline M^t&N^t& A&B\\
-i\overline N^t&-M^t&-B^*&-A^t
\end{array}\right)$$
} where $\a\in i\bbbr,$ $M,N\in M_{1\times n}(\bbbc),$ $ A,B\in M_{n\times n}(\bbbc),$  $B^t=B,$ and $A^*=-A.$ {The Cartan--Killing form $\kappa$  of $\frak{c}(n)$ is a real nonzero scalar multiple of its supertrace form. All derivations of $\frak{c}(n)$ are inner because $\kappa$ is nondegenerate, in particular $H^2(\frak{c}(n))=\{0\}.$  The compact Lie algebra  $\frak{sp}(n-1)$ has a negative definite Killing form and  the trace form on $\frak{sp}(n-1)$ is a positive real scalar multiple  of the Killing form.}  Therefore $\kappa$ restricted to $\bbbr\sub \frak{c}(n)_0$ is negative definite while the restriction of $\kappa$ to $\frak{sp}(n-1)\sub \frak{c}(n)_0$ is positive definite. So it is easy to {find}    $x_1\in \bbbr$ and $x_2\in \frak{sp}(n-1)$ such that
{
\begin{equation}\label{equalzero2}
\kappa(x_1,x_2)=0\andd \kappa(x_1,x_1)=-\kappa(x_2,x_2)\neq 0;
\end{equation}
 in particular, $\kappa(x_1+x_2,x_1+x_2)= 0.$}
\medskip

 \noindent{{\bf III}. $\pmb{\frak{pq}(n)}$   $\pmb{(n>2)}:$} Suppose $n>2.$ Set $$\frak{q}(n)=\left\{\left(\begin{array}{cc}a&(1-i)b\\
(1-i)b&a\end{array}\right)\vline~~ a,b\in\fu(n;\bbbc),{\rm tr}(b)=0\right\}$$ and  take $$\frak{pq}(n):=\frak{q}(n)/\bbbr i{\bf 1}\simeq\left\{\left(\begin{array}{cc}a&(1-i)b\\
(1-i)b&a\end{array}\right)\vline~~ a,b\in\fsu(n;\bbbc)\right\}.$$
The real Lie superalgebra $\frak{q}(n)$ is a real form of the  Lie superalgebra $$\widetilde Q(n):=\left\{\left(\begin{array}{cc}
a&b\\b&a
\end{array}\right)\in \frak{sl}(n+1,n+1)
\vline~~ {\rm tr}(b)=0\right\}$$ and $\frak{pq}(n)$ is a real form of $Q(n):=\widetilde Q(n)/\bbbc{\bf 1}.$   The even part {$\frak{pq}(n)_0\simeq \frak{su}(n,\bbbc)$} is a simple Lie algebra acting irreducibly on {$\frak{pq}(n)_1\simeq \frak{su}(n,\bbbc),$ by the adjoint representation.} The restriction of {$D:={\rm ad}\left(\begin{array}{cc}
0&{\bf 1}_n\\i{\bf 1}_n&0
\end{array}\right)$ to $\frak{pq}(n),$ which is an outer derivation of $\frak{pq}(n)$ vanishing on $\frak{pq}(n)_0,$}  forms  a  basis for $H^2(\frak{pq}(n))$ regarding Remark \ref{rem2}. Also the bilinear form
\begin{align*}
\kappa\colon &\frak{pq}(n)\times  \frak{pq}(n)\longrightarrow \bbbr\\
&\left( \left(\begin{array}{cc}a&(1-i)b\\
(1-i)b&a\end{array}\right),\left(\begin{array}{cc}a'&(1-i)b'\\
(1-i)b'&a'\end{array}\right)\right)\mapsto {\rm tr}(ab'+a'b)
\end{align*}
is an odd invariant nondegenerate supersymmetric bilinear form on $\frak{pq}(n).$  It is easily verified  that the  bilinear form   {{$\kappa_D\colon \fk\times \fk\longrightarrow \bbbr,$  $(x,y)\mapsto \kappa(D(x),y),$} is  definite on $\fk_1\times \fk_1.$}

\begin{Thm}[{\cite[Thm. 2.3]{AN}}] \label{S-K} Each simple compact Lie superalgebra $\fk$
is either a simple compact Lie algebra or isomorphic to one  of the Lie superalgebras
  \begin{equation}
    \label{eq:comp-super}
\fsu(n|m;\bbbc), n>m,\quad \mbox{ and }
\quad \frak{psu}(n|n;\bbbc), \quad\frak{pq}(n), \quad \frak{c}(n), n\geq 2.
  \end{equation}
\end{Thm}

Note that we have seen in (I)-(III) above that any of these
Lie superalgebras $\fk$ carries a nondegenerate supersymmetric invariant homogeneous bilinear form.

\begin{rem}\label{rem4}{\rm
(i) It is easily verified that, if a real vector space $V$ equipped with a bilinear form $\fm$ has an orthogonal  decomposition $V=V_1\op V_2$ of  nonzero subspaces such that $\fm$ is negative definite on $V_1$ and positive definite on $V_2,$ then $\{x\in V|(x,x)=0\}$ spans the vector space $V.$ Now if $\fk$ is one of the real Lie superalgebras $\fsu(n|m;\bbbc)$ $(n>m),$ $\frak{psu}(n|n;\bbbc)$ $(n\geq 2)$ or $\frak{c}(n)\;(n\geq 2)$ and $\kappa$ is the bilinear form introduced in (I) and  (II), then $$\{x\in\fk_0\mid \kappa(x,x)=0\}$$ spans $\fk_0.$

(ii) Let  $\fk$ be one of the simple compact Lie superalgebras in
\eqref{eq:comp-super} and $\kappa$ be the nondegenerate supersymmetric
invariant bilinear form on $\fk$ introduced in (I)--(III).
Since the real Lie superalgebra  $\fk$ is a real form of a simple Lie superalgebra, it is absolutely simple and so $\cent_\bbbr(\fk)=\bbbr{\rm id}.$ In particular, by Proposition~\ref{cor-base}, up to scalar multiple, $\kappa$ is the unique nonzero supersymmetric invariant bilinear form on~$\fk.$

(iii) Suppose that $\fk$ is one of the simple compact Lie superalgebras
$\fsu(n|m;\bbbc)$ with $n>m$ or  $\frak{c}(n)$ with $n\geq 2.$ We know
from (I) and (II) that $\fk_0$ is a reductive Lie algebra  and  the center
$Z(\fk_0)$  of $\fk_0$ is isomorphic to $\bbbr.$ A  direct calculation shows that the  bilinear form $b_\frak{z}\colon \fk_1\times\fk_1\longrightarrow \bbbr$ mapping $(x,y)$ to the projection of $[x,y]$ on the center component of $[x,y]\in\fk_0,$ with respect to the decomposition of $\fk_0$ stated in (I) and (II), is a definite form. This in particular implies that $\sc(\fk)$ is pointed, (Lemma~\ref{point-crit};
see also \cite[Thm. 5.4 \& Lem. 3.4]{AN}).

(iv) Suppose that $n\geq 2$ is a positive integer. Then for
\[ B_j:=\diag(0,\ldots,0,\underbrace{1}_{j{\rm th}},0,\ldots,0)\quad\quad\hbox{for $1\leq j\leq n$}\]
  and  $$X_j:=\left(
\begin{array}{cc}
0&B_j\\
iB_j&0
\end{array}
\right)\in \frak{su}(n;\bbbc)_1\quad\quad\quad \hbox{for $1\leq j\leq n$},$$  we have $\sum_{j=1}^n[X_j,X_j]\in i\bbbr {\bf 1}.$ Also
for $b_j:=i(B_j-B_{j+1})$ $(1\leq j\leq n-1)$ and   $b_n:=i(B_n-B_1),$ we have
$$Y_j:=\left(
\begin{array}{cc}
0&(1-i)b_j\\
(1-i)b_j^*&0
\end{array}
\right)\in \frak{pq}(n;\bbbc)_1\quad\quad\quad (1\leq j\leq n)$$ and  $\sum_{j=1}^n[Y_j,Y_j]\in i\bbbr{\bf 1}.$ This implies that,  for $\fk=\frak{psu}(n|n;\bbbc)$ $(n\geq 2)$ or $\fk=\frak{pq}(n)$ $(n>2),$ $\sc(\fk)$ is not pointed.
}
\end{rem}

\begin{Thm}\label{unitary-rep}
Suppose that $\fk$ is  a simple compact Lie superalgebra with $\fk_1\neq \{0\}$ and $s$ is a positive integer. Then {$\Lam_s^+(\bbbr)\ot \fk$} lies in  the kernel of each unitary representations of $\fg=\Lam_s(\bbbr)\ot \fk;$ in particular  unitary representations {of  $\fg$  are in one-to-one  correspondence with}    unitary representations of $\fk.$
\end{Thm}

\begin{prf} We first note that  $\fg= ( \Lam_s^+(\bbbr)\ot \fk)\oplus \fk$
and that $ \Lam_s^+(\bbbr)\ot \fk$ is an ideal of $\fg.$  So each
unitary representation $\pi$ {of}  $\fk$ can be extended to a unitary representation {of}  $\fg$ with
$\Lam_s^+(\bbbr)\ot \fk\sub{\rm ker}(\pi).$

{Next recall that $\fk_0$ is a reductive Lie algebra and set $\fh:=[\fk_0,\fk_0].$ Suppose that $\pi$ is a unitary representation of $\fg.$ Then for a nonnegative integer $t$ and elements  $i_1,\ldots,i_{2t+1}\in\{1,\ldots,s\},$ thanks to Lemma \ref{necessary}(i), we have
\begin{equation}\label{in-ker3}
\ep_{i_1}\cdots\ep_{i_{2t+1}}\ot\fk_0\sub {\rm ker}(\pi).
\end{equation} Therefore,  for $x,y\in \fh$ and  $i_1,\ldots,i_{2t+2}\in\{1,\ldots,s\},$ $$ \ep_{i_1}\cdots\ep_{i_{2t+2}}\ot[x,y]=[\ep_{i_1}\ep_{i_2}\cdots\ep_{i_{2t+1}}\ot x, \ep_{i_{2t+2}}\ot y]\in {\rm ker}(\pi).$$
So as $\fh$ is perfect, we get using (\ref{in-ker3}) that $${\Lam^+_s(\bbbr)\ot \fh\sub{\rm ker}(\pi)}.$$ {Fixing} $x\in \fh$ and $y\in \fk_1$ with $[x,y]\neq 0,$ we have  for each $i_1,\ldots,i_t\in \{1,\ldots,s\},$
$$0\neq\ep_{i_1}\cdots\ep_{i_{t}}\ot[y,x]=[ 1\ot y, \ep_{i_1}\cdots\ep_{i_{t}}\ot x]\in {\rm ker}(\pi)\cap (\Lam_s^t(\bbbr)\ot \fk_1).$$ So it follows that $${\Lam^+_s(\bbbr)\ot \fk_1\sub{\rm ker}(\pi)}$$  because  $\fk_1$ is an irreducible $\fk_0$-module. We finally note {that,} as $\fk$ is simple, $\fk_0=[\fk_1,\fk_1]$ and so
\begin{align*}
\Lam^+_s(\bbbr)\ot \fk_0=\Lam^+_s(\bbbr)\ot [\fk_1,\fk_1]=[1\ot \fk_1,\Lam^+_s(\bbbr)\ot \fk_1]\sub[\fg,{\rm ker}(\pi)]
\sub {\rm ker}(\pi).
\end{align*}
This completes the proof.}
\end{prf}

One knows  from \cite[Lem.'s~3.2 \&~3.4(v)]{AN} that faithful finite dimensional
 unitary representations do not exist  for
$\frak{psu}(n|n;\bbbc)$ and $\frak{pq}(n)$ ($n>2$)
while $\frak{su}(n;\bbbc)$ and $\frak{q}(n)$ ($n>2$) which are respectively
universal central extensions of $\frak{psu}(n|n;\bbbc)$ and $\frak{pq}(n)$ ($n>2$)
have finite dimensional  faithful unitary representations.
We also know form  Theorem \ref{unitary-rep} that there is no faithful
unitary representation for  $\fg=\Lam_s(\bbbr)\ot \fk$ but we are interested in
faithful unitary representations of  central extensions of $\fg.$

In what follows  recalling (\ref{equalzero1}) and (\ref{equalzero2}), we will  see that
for each  central  extension $\widehat\fg$ of $\fg,$ the ideal
${\rm urad}(\widehat\fg)$ is non-zero.
In particular, $\widehat\fg$ does not have faithful unitary representations
(Lemma \ref{necessary}).

\begin{Thm}\label{urad2}
Let $s$ be a positive integer and { $\fk$ be a compact simple} Lie superalgebra with $\fk_1\neq\{0\}.$  Suppose  $\widehat\fg=\g\op M$
is a  perfect  central extension of $\fg=\Lam_s(\bbbr)\ot \fk$ such that    ${\rm urad}(\widehat\fg)$ is a proper ideal, then  $${\rm urad}(\widehat\fg)=(\Lam^+_s(\bbbr)\ot \fk)\op (M\cap {\rm urad}(\widehat\fg)).$$
In particular, either all unitary representations of $\hat\g$ are trivial,
or they factor through a central extension of $\fk$.
\end{Thm}
\begin{prf}
Suppose that  $\omega$ is the $2$-cocycle corresponding to $\widehat\fg.$
We split the proof into three steps:

\smallskip

\noindent{\bf Step 1.} $\Lam^+_s(\bbbr)\ot \fk\sub {\rm urad}(\hat\fg):$ To show this, we use a type--by--type approach:

\smallskip

\noindent\underline{${\fk=\frak{pq}(n)}$:}  Consider the  outer derivation $D$ of $\fk$ introduced in (III). By Remark \ref{rem4}(ii) and Theorem \ref{cor1}, there {exist} a linear map $f$ on $\Lam_s(\bbbr)$ and a Hochschild map $F$ on $\Lam_s(\bbbr)\times\Lam_s(\bbbr)$ such that for $a,b\in\Lam_s(\bbbr)$ and $x,y\in\fk,$  $$\omega(a\ot x,b\ot y)=(-1)^{|x||b|}\Big(F(a,b)\kappa(x,y)+f(ab)\kappa(Dx,y)\Big).$$

\noindent $\bullet$ Stage 1. $\Lam^{\rm odd}_s(\bbbr)\ot \fk_0\sub {\rm urad}(\widehat\fg):$
As $\kappa$ is  odd and $D\mid_{\fk_0}=0,$ we have
\begin{equation}\label{even-part}
[a\ot x,b\ot y]_\omega=ab\ot[x,y]\quad \hbox{for $x,y\in \fk_0,$  $a,b\in \Lam_s(\bbbr).$}
\end{equation}
 This  shows that, for $a\in \Lam^{2m+1}(\bbbr)$ $(m\geq 1)$ and $x\in\fk_0,$ $a\ot x$ is an odd element of $\widehat\fg$ with square zero  and so by Lemma \ref{necessary}(i), $$
 \Lam^{\rm odd}_s(\bbbr)\ot \fk_0\sub {\rm urad}(\widehat\fg).
 $$
 \noindent $\bullet$  Stage 2. $\Lam^+_s(\bbbr)\ot \fk_0\sub{\rm urad}(\widehat\fg):$
 Since $\fk_0$ is simple, Stage 1 together with (\ref{even-part}) implies that
\begin{align*}
\Lam^{\rm even}_s(\bbbr)\ot \fk_0
&= \sum_{m=1}^\infty\Lam^{2m-1}_s(\bbbr)\Lam^{1}_s(\bbbr)\ot [\fk_0,\fk_0]\\
&=\sum_{m=1}^\infty[\Lam^{2m-1}_s(\bbbr)\ot\fk_0,\Lam^{1}_s(\bbbr)\ot\fk_0]_\omega
 \sub {\rm urad}(\widehat\fg).
 \end{align*}
 Therefore, we get the result using Stage 1.

 \noindent $\bullet$ Stage 3. $\Lam_s^+(\bbbr)\ot \fk\sub {\rm urad}(\widehat\fg):$ For $x\in \fk_0,$ $y\in \fk_1$ and $a\in \Lam_s^+(\bbbr),$ by (\ref{Hochschild}) and Stage 2, we have
\begin{fleqn}\begin{align*}
a\ot[x,y]=a\ot[x,y]+F(a,1)\kappa(x,y)+f(a)\kappa(Dx,y)=&[a\ot x,1\ot y]_\omega\\
\in& [{\rm urad}(\widehat\fg),\widehat{\fg}]_\omega
\sub {\rm urad}(\widehat\fg).
\end{align*}
\end{fleqn}
So choosing $x\in \fk_0$ and $y\in \fk_1$ with $u:=[x,y]\neq 0$ ($\fk_1$ is an irreducible $\fk_0$-module), we have $\Lam_s^+(\bbbr)\ot u\sub {\rm urad}(\widehat\fg).$ Now for $x_1,\ldots,x_k\in \fk_0,$ we have
\begin{fleqn}
\begin{align*}
\Lam^+(\bbbr)\ot[x_k,[\ldots,[x_1,u]\ldots]]=&[ 1\ot x_k,[\ldots,[1\ot x_1,\Lam_s^+(\bbbr)\ot u]_\omega\ldots]_\omega]_\omega\nonumber\\
\sub& [\widehat{\fg},{\rm urad}(\widehat\fg)]_\omega\sub {\rm urad}(\widehat\fg).
\end{align*}
\end{fleqn}
As $\fk_1$ is irreducible, this implies that $\Lam_s^+(\bbbr)\ot \fk_1\sub {\rm urad}(\widehat\fg)$ and so  Stage 2 completes the proof.

\noindent \underline{${\fk=\frak{psu}(n|n;\bbbc)}:$} Consider the  outer derivation $D$ of $\fk$ introduced in (I). By Remark \ref{rem4}(ii) and Theorem \ref{cor1}, there are a linear map $f$ on $\Lam_s(\bbbr)$ and a Hochschild map $F$ on $\Lam_s(\bbbr)\times\Lam_s(\bbbr)$ such that for $a,b\in\Lam_s(\bbbr)$ and $x,y\in\fk,$  $$\omega(a\ot x,b\ot y)=(-1)^{|x||b|}\Big(F(a,b)\kappa(x,y)+f(ab)\kappa(Dx,y)\Big).$$

\noindent $\bullet$ Stage 1. $\Lam^{\rm odd}_s(\bbbr)\ot \fk_0\sub {\rm urad}(\widehat\fg):$
Suppose $a\in \Lam_s(\bbbr)$ and   $x\in \fk_0$ with  $\kappa(x,x)=0.$ {Then,} as $Dx=0,$ we have
\begin{equation}\label{even-part1}
[a\ot x,a\ot x]_\omega=a^2\ot[x,x]+F(a,a)\kappa(x,x)+f(a^2)\kappa(Dx,x)=0.
\end{equation}
So by Lemma \ref{necessary}(i),  we get $$\{a\ot x\mid a\in \Lam^{\rm odd}_s(\bbbr),x\in\fk_0,\kappa(x,x)=0 \}\sub {\rm urad}(\widehat\fg).$$ Therefore, by Remark \ref{rem4}(i), we have
$$
 \Lam^{\rm odd}_s(\bbbr)\ot \fk_0\sub {\rm urad}(\widehat\fg).
 $$

\noindent $\bullet$ Stage 2. $\Lam_s^+(\bbbr)\ot\fk_1\sub {\rm urad}(\widehat\fg):$
  Since $\kappa$ is even and $D\mid_{\fk_0}=0,$ we have
\begin{equation*}
[a\ot x,b\ot y]_\omega=ab\ot[x,y]\quad\quad\hbox{ for $x\in\fk_0,y\in\fk_1,$ $a,b\in \Lam_s(\bbbr).$}
\end{equation*} So
Stage 1 implies that for $x,x_1,\ldots,x_k\in \fk_0$ and $y\in\fk_1,$
\begin{multline*}
\Lam_s^+(\bbbr)\ot[x_k,[\ldots,[x_1,[x,y]]\ldots]]\\
\begin{array}{l}
=[1\ot x_k,[\ldots,[1\ot x_1,[\Lam_s^{\rm odd}(\bbbr)\ot x,\Lam_s(\bbbr)\ot y]_\omega]_\omega\ldots]_\omega]_\omega\\
\sub [\widehat{\fg},{\rm urad}(\widehat\fg)]_\omega\sub {\rm urad}(\widehat\fg).
\end{array}
\end{multline*}
 Now we are done  as $\fk_1$ is an irreducible $\fk_0$-module.

 \noindent $\bullet$ Stage 3. $\Lam_s^+(\bbbr)\ot \fk\sub {\rm urad}(\widehat\fg):$
Recalling $x_*$ and $y_*$ from (\ref{special}) and using (\ref{f-use}), we have
\begin{equation}\label{star}
\Lam^{\rm even}_s(\bbbr)\ot[x_*,y_*]=[\Lam^{\rm odd}\ot x_*,\Lam^{\rm odd}\ot y_*]_\omega\sub {\rm urad}(\widehat\fg).
\end{equation} We recall that ${\frak psu}(n|n;\bbbc)_0=\fk_0^1\op\fk_0^2$ with $\fk_0^1,\fk_0^2\simeq\frak{su}(n;\bbbc)$ and that $[x_*,y_*]=u+v,$ $0\neq u\in\fk_0^1$ and $0\neq v\in \fk_0^2.$ Now for $x_1,\ldots,x_k,x\in\fk_0^1,$ using (\ref{star}), we have
\begin{multline*}
\Lam_s^{\rm even}(\bbbr)\ot[x_k,[\ldots,[x_1,[x,u]]\ldots]]\\
\begin{array}{l}=[1\ot x_k,[\ldots,[1\ot x_1,[1\ot  x,\Lam_s^{\rm even}(\bbbr)\ot (u+v)]_\omega\ldots]_\omega]_\omega\\
\sub [\widehat{\fg},{\rm urad}(\widehat\fg)]_\omega
\sub {\rm urad}(\widehat\fg).
\end{array}
\end{multline*}
Also, for $y_1,\ldots,y_k,y\in\fk_0^2,$ we have
\begin{multline*}
\Lam_s^{\rm even}(\bbbr)\ot[y_k,[\ldots,[y_1,[y,v]]\ldots]]\\
\begin{array}{l}
=[1\ot y_k,[\ldots,[1\ot y_1,[ 1\ot y,\Lam_s^{\rm even}(\bbbr)\ot (u+v)]_\omega\ldots]_\omega]_\omega\\
\sub [\widehat{\fg},{\rm urad}(\widehat\fg)]_\omega
\sub {\rm urad}(\widehat\fg).
\end{array}
\end{multline*} Since $\fk_0^1$ and $\fk_0^2$ are simple, these imply that $\Lam^+_s(\bbbr)\ot \fk_0\sub {\rm urad}(\widehat\fg).$ This together with Stages 1--2  gives that $\Lam_s^+(\bbbr)\ot \fk\sub {\rm urad}(\widehat\fg).$

\noindent \underline{${\fk=\frak{su}(m|n;\bbbc),m\neq n}$:} In this case, $\fk_0=\fk_0^1\op\fk_0^2\op\bbbr i\i$ where $i\i$ is a central element and $\fk_0^1$ as well as $\fk_0^2$ are simple ideals of $\fk_0.$ Since   $H^2(\fk)=\{0\}$ by   (I), using
Remark~\ref{rem4}(ii) and  Theorem \ref{cor1}, there exists  a  Hochschild map $F$ on $\Lam_s(\bbbr)\times \Lam_s(\bbbr) $ such that for $a,b\in \Lam_s(\bbbr)$ and $x,y\in \fk,$  $$\omega(ax,by)=(-1)^{|x||b|}F(a,b)\kappa(x,y).$$

\noindent $\bullet$ Stage 1. $ (\Lam^{\rm odd}_s(\bbbr)\ot \fk_0)\op(\Lam^+_s(\bbbr)\ot\fk_1)\sub {\rm urad}(\widehat\fg):$ For $a\in \Lam_s(\bbbr)$ and  $ x\in\fk_0$  with $\kappa(x,x)=0,$ we have $$[a\ot x,a\ot x]_\omega=a^2\ot[x,x]+F(a,a)\kappa(x,x)=0.$$
A modification of the argument  in Stages 1--2 of the previous case implies that
 $$
 (\Lam^{\rm odd}_s(\bbbr)\ot \fk_0)\op (\Lam^+_s(\bbbr)\ot\fk_1)\sub {\rm urad}(\widehat\fg).
 $$

 \noindent $\bullet$ Stage 2. $ \Lam^{\rm even}_s(\bbbr)\ot (\fk_0^1\op\fk_0^2)\sub {\rm urad}(\widehat\fg):$ Choose $z_*,u,v$ and $z$ as in (\ref{z*}) and (\ref{uv}) and
use a slight modification of the argument in Stage 3 of the previous case.

 \smallskip

 \noindent $\bullet$ Stage 3. $\Lam^+_s(\bbbr)\ot\fk\sub {\rm urad}(\widehat\fg):$ Choose $x\in\fk^1_0$ with $\kappa(x,x)\neq 0,$  then for $b\in\Lam^+_s(\bbbr)$ and $a\in \Lam_s(\bbbr),$  by Stages 1--2, we have
\[
 \kappa(x,x)F(a,b)=[a\ot x,b\ot x]_\omega\in [\hat\fg,{\rm urad}(\widehat\fg)]_\omega\sub{\rm urad}(\widehat\fg).
 \]
It means that
 \begin{equation*}\label{cent-comp}
F(a,b)\in{\rm urad}(\widehat\fg)\quad\quad\hbox{for $a\in\Lam_s(\bbbr)$  and $b\in \Lam_s^+(\bbbr).$}
 \end{equation*}
This together with Stage 1 implies that for $y\in\fk_1,$ we have
\begin{eqnarray}
  \label{fab}
\Lam_s^{\rm even}(\bbbr)\ot[y,y]&\sub&[\Lam_s^{\rm odd}(\bbbr)\ot y,\Lam_s^{\rm odd}(\bbbr)\ot y]_\omega+F(\Lam_s^{\rm odd}(\bbbr),\Lam_s^{\rm odd}(\bbbr))\kappa(y,y)\nonumber\\
&\sub& {\rm urad}(\widehat\fg).
\end{eqnarray}
Fix a nonzero $y\in\fk_1.$ As we mentioned in Remark \ref{rem4}(iii), $[y,y]$ has a nonzero component in $\bbbr i\i.$ So by (\ref{fab}) and Stage 2, we get that $\Lam_s^{\rm even}(\bbbr)\ot\bbbr i\i\sub {\rm urad}(\widehat\fg).$ This together with Stages 1--2,  completes the proof.
\smallskip

\noindent \underline{${\fk=\frak{c}(n), n\geq 2}:$}
We recall from (II) that  $\fk_0=\frak{sp}(n-1)\op \bbbr.$ Since $\fk$ is simple, $\fk_0=[\fk_1,\fk_1]$ and so there are odd elements $x,y\in\fk$ such that $[x,y]$ has a nonzero component in $\frak{sp}(n-1).$ So using the same argument as in the previous case, we get  the result.

\medskip

\noindent{\bf Step 2.} ${\rm urad}(\widehat\fg)\cap ((\bbbr\ot \fk)\op M)
\subeq M:$
We first suppose $\fk=\frak{pq}(n)$ and consider $f$ as above. If  $(1\ot x_0)+m_0\in{\rm urad}(\widehat\fg)$ for  some
$m_0\in M$ and  nonzero $x_0\in \fk,$ then we have for $y \in \fk$
\begin{equation}
  \label{eq:e1}
1\ot [x_0,y]+f(1)\kappa(Dx_0,y)
=[1\ot x_0,1\ot y]_\omega =[1\ot x_0+m_0,1\ot y]_\omega
\in {\rm urad}(\widehat\fg).
\end{equation}
But $\fk$ is simple, so this implies that, for each $x\in \fk,$ there is $r_x\in M$ with
\begin{equation}\label{final2}
(1\ot x)+r_x\in {\rm urad}(\widehat\fg).
\end{equation}
 Since $\widehat\fg$ is perfect,  for $x\in \fk$ and $m\in M,$ there are $x_1,\ldots,x_n,y_1,\ldots,y_n\in \fk$ and homogeneous elements $a_1,\ldots,a_n,b_1,\ldots,b_n\in \Lam_s(\bbbr)$ with $a_i=1$ if $a_i\in \Lam_s^0(\bbbr),$ such that $x+m=\sum_{i=1}^n[a_i\ot x_i,b_i\ot y_i]_\omega.$ So as $\Lam_s^+(\bbbr)\ot \fk\sub{\rm urad}(\hat\fg),$ (\ref{final2}) implies that
 \begin{align*}
 x+m
 =&\sum_{\substack{i=1\\
 a_i\neq 1}
 }^n[a_i\ot x_i,b_i\ot y_i]_\omega
 +\sum_{\substack{i=1\\
 a_i= 1}
 }^n[a_i\ot x_i,b_i\ot y_i]_\omega\\
 =&\sum_{\substack{i=1\\
 a_i\neq 1}
 }^n[a_i\ot x_i,b_i\ot y_i]_\omega
 +\sum_{\substack{i=1\\
 a_i= 1}
 }^n[(a_i\ot x_i)+r_{x_i},b_i\ot y_i]_\omega \in {\rm urad}(\widehat\fg)
 \end{align*}
 which is a contradiction as ${\rm urad}(\widehat\fg)$ is proper. So we are done  in this case. Repeating this argument, one gets the result  for $\fk=\frak{psu}(n|n;\bbbc)$ as well.

Next assume $\fk$ is one of the remaining types.
If  $(1\ot x_0)+m_0\in{\rm urad}(\widehat\fg)$ for  some $m_0\in M$ and  nonzero $x_0\in \fk,$ then for $y\in \fk,$ we have
\begin{align*}
1\ot [x_0,y]=[1\ot x_0,1\ot y]_\omega=[1\ot x_0+m_0,1\ot y]_\omega\in {\rm urad}(\widehat\fg).
\end{align*}
But $\fk$ is simple, so it follows  that  $1\ot  \fk\sub {\rm urad}(\widehat\fg).$ This implies that  $\Lam_s(\bbbr)\ot \fk\sub {\rm urad}(\widehat\fg)$ because $\Lam^+_s(\bbbr)\ot \fk\sub {\rm urad}(\widehat\fg).$ Since $\widehat\fg$ is perfect,   we get that
\[ \widehat\fg=[\Lam_s(\bbbr)\ot \fk,\Lam_s(\bbbr)\ot \fk]_\omega
\sub {\rm urad}(\widehat\fg).\]
This  is  a contradiction and so we are done.
\medskip

\noindent{\bf Step 3.} ${\rm urad}(\widehat\fg)=(\Lam^+_s(\bbbr)\ot \fk)\op (M\cap {\rm urad}(\widehat\fg)):$ By Step 1,  $\Lam^+_s(\bbbr)\ot \fk\sub{\rm urad}(\hat\fg).$ Therefore, we have
\[ {\rm urad}(\widehat\fg)=(\Lam_s^+(\bbbr)\ot \fk)\op
({\rm urad}(\widehat\fg)\cap ((\bbbr\ot \fk)\op M))\] and so we  are done by Step 2.
\end{prf}

\appendix
\section{Unitary representations of semidirect products}
\label{sec:3}

In this section, we describe the classification of irreducible
unitary representations of Lie supergroups which are
semidirect products $\cN \rtimes K$ of a finite dimensional
Clifford--Lie supergroup $\cN$ and a compact Lie group $K$
(cf.\ Definition~\ref{def:lie-subgrp}).
This classification has been obtained in \cite{CCTV}, but since it is also
central in our context, we recall this result is some detail,
and this requires a number of ingredients, such as representations
of Clifford algebras and Mackey's little group theory for Lie
supergroups.

\subsection{Clifford algebras}
Let $V$ be a  finite dimensional real vector space and  $\mu$ a  symmetric bilinear form on  $V.$  We write $\cC = \cC(V,\mu)$ for the corresponding Clifford algebra and
$\iota \: V \to \cC$ for the structure map satisfying
\[ \iota(v)^2 = \mu(v,v) \one \quad \mbox{ for } \quad v \in V.\]
As $\iota$ is injective, we shall identify $V$ with the subset $\iota(V) \subeq \cC$.
We consider the $\bbbz_2$-grading  $\cc=\cc_0\op\cc_1$ on $\cc$ induced  from the natural $\bbbz$-grading on $\cc$ and recall  the {\it parity operator}
\[ \Pi_\cC:\cc\longrightarrow \cc, \quad x_0 + x_1 \mapsto x_0 - x_1.\]
The unit group $\cC^\times$ acts on $\cC$ by the automorphisms
\begin{equation}
  \label{eq:cliffgrp-act}
 \alpha_g(c) := \Pi_\cC(g) c g^{-1}, \quad x \in \cC^\times, c \in \cC.
\end{equation}
The stabilizer of $V$ is the {\it Clifford group}
\[ \Gamma(V) := \{ g \in \cC^\times \:
\alpha_g(V) = V\}.\]
%If $\mu$ is non-degenerate, then it is adapted to the grading in the sense that
%$\Gamma(V)=(\Gamma(V)\cap\cc_0)\cup(\Gamma(V)\cap\cc_1)$ and
The map
\[ \a:\Gamma(V)\longrightarrow \OO(V),\quad\quad
x\mapsto \a_x\res_V \]
defines a group homomorphism.
Let $x \mapsto x^T$ denote the unique antiautomorphism of $\cC$
that coincides with the identity on $V$.
Then $x^Tx = N(x)\one$ for $x\in \Gamma(V)$ and some $N(x) \in \R^\times$,
which defines a group homomorphism
$N:\Gamma(V)\longrightarrow \bbbr^\times$.
Its kernel is the  {\it Pin group}  $\Pin(V)=\Pin(V,\mu):={\rm ker}(N)$.
%It inherits the decomposition
%\[ \Pin(V)=(\Pin(V)\cap\cc_0)\cup (\Pin(V)\cap\cc_1).\]
If $\mu$ is positive definite, then each element of $\Pin(V)$ is a
product $v_1\cdots v_n$ of unit length vectors $v_j \in V$ and
%$\a$ in the sequence (\ref{sequence}) is surjective and
we obtain the following short exact sequence
\begin{equation}\label{sequence}{\bf 1} \longrightarrow \{\pm 1\}\stackrel{}{\longrightarrow} \Pin(V)\stackrel{\a}{\longrightarrow} \OO(V){\longrightarrow} {\bf 1}.
\end{equation}
%$${\bf 1}\longrightarrow \{\pm 1\}\stackrel{N}{\longrightarrow} {\rm Spin}(V)\stackrel{\a}{\longrightarrow} S\OO(V).$$

The universal property of Clifford algebras implies that each element $T\in \OO(V)$
induces an automorphism $\nu_T$ of $\cc$ with $\nu_T\res_V = T$.
This defines a natural action $\nu: \OO(V)\longrightarrow {\rm Aut}(\cc)$.
%The corresponding action $\Ad^p := \nu\circ\a$
%of $\Pin(V)$ on $\cc$ is the {\it parity-adjoint action}
%given by
%\begin{equation}\label{ad-parity}
%\begin{array}{rl}
%\Ad^p_x(y) = {(-1)^{|x|}xyx^{-1}} \quad \mbox{ for }\quad
%x \in \Pin(V), y \in \cC.
%\end{array}
%\end{equation}
% KH: it seems that $\Ad^p is never used below.

\begin{deft} {\rm A {\it selfadjoint representation} $(\pi, \hH)$ of a
Clifford algebra $\cC = \cC(V,\mu)$ on a Hilbert superspace $\hH$
is an algebra homomorphism $\pi \: \cC \to \End(\hH)$ for which all operators
$\pi(v)$, $v \in V$, are odd and symmetric. Therefore selfadjoint representations
are in {one-to-one} correspondence with linear maps \break
$\pi \: V \to \Herm_1(\hH)_{\rm odd}$ satisfying
\begin{equation}
  \label{eq:square-rel}
 \pi(v)^2 = \mu(v,v)\one \quad \mbox{ for } \quad v \in V.
\end{equation}}
\end{deft}

By \cite[Lem. 11]{CCTV}, we have:

\begin{Pro}\label{unique} If $\mu$ is positive definite,
then there exists a finite dimensional irreducible selfadjoint representation
$(\tau,\hN)$ of $\cC(V,\mu)$ which is unique
if $\dim V$ is odd and unique up to parity reversal if $\dim V$ is even.
Any other selfadjoint representation
$(\pi,\hH)$ is of the form $\pi = \one \otimes \tau$, where
$\hH = \hM \otimes \hN$ is a tensor product of Hilbert superspaces.
If $\dim V$ is odd, then the multiplicity space $\hM$ can be chosen purely even.
\end{Pro}

\begin{deft}{\rm
If $G$ is a group and $\hH$ a Hilbert superspace, then
a homomorphism $\pi \: G \to \Aut(\hH)$ is called a {\it graded
unitary representation} with respect to the subgroup  {$G_1:= \pi^{-1}(\Aut(\cH)_{\rm even})$  if $G_1$} is a proper subgroup.
It is called {\it even} if $\pi(G) \subeq \Aut(\hH)_{\rm even}$, i.e.,
if $G_1 = G$.}
\end{deft}

\begin{rem}\label{rem:super-ext}
{\rm If $G_1 \subeq G$ is a subgroup of index two (hence normal),
any unitary representation
$(\pi, \hH_0)$ of $G_1$ on a (purely even) Hilbert space $\hH_0$
admits, up to equivalence, a unique extension to a
graded representation $(\hat\pi, \hH)$ of $G$ on a
Hilbert superspace $\hH = \hH_0 \oplus \hH_1$. It is equivalent
to $\Ind_{G_1}^G(\pi)$ with its natural $2$-grading.
%For a fixed element
%$g \in G \setminus G_1$ and $\tau(h) := g^{-1}hg$, it is equivalent
%to the unitary representation on $\hH_0 \oplus \hH_0$ given by
%\[ \hat\pi(g) = \pmat{0 & \pi(g^2)\\ \one & 0}, \quad
%\hat\pi(h) = \pmat{ \pi(h) &  0 \\ 0 & \pi(\tau(h))}
%\quad \mbox{ for }\quad h \in G_1.\]
}\end{rem}

In the following proposition we formulate the outcome of Lemmas 13/14
in \cite{CCTV} in terms of a central extension of the orthogonal group.

\begin{Pro} \label{prop:ohat}
Suppose that $\mu$ is positive definite and that
$(\tau,\hN)$ is an irreducible selfadjoint representation of
$\cC(V,\mu)$. Then there exists a
central extension $q_O \: \hat\OO(V,\mu) \to \OO(V,\mu), g \mapsto \oline g$
of $\OO(V,\mu)$ by the two-element
group $\{\pm 1\}$ such that:
\begin{itemize}
\item[\rm(i)]  The subgroup $q_O^{-1}(\SO(V,\mu))$ is equivalent, as a central
extension, to \break $\Spin(V) := \alpha^{-1}(\SO(V))\subeq \Pin(V)$.
\item[\rm(ii)] There exists a graded unitary representation
$\hat\tau \:  \hat\OO(V,\mu) \to \Aut(\hN)$
extending the even representation $\tau\res_{\Spin(V)}$
of $\Spin(V,\mu)$ such that
\[  \hat\tau(g)\tau(c) \hat\tau(g)^{-1} = \tau(\nu_{\oline g}(c))
\quad \mbox{ for } \quad
g \in \hat\OO(V,\mu), c \in \cC(V,\mu).\]
If $\dim V$ is odd, then there exists an even representation
$\hat\tau$ with this property.
\end{itemize}
\end{Pro}

\subsection{Clifford--Lie superalgebras}

We now explain how selfadjoint representations of Clifford algebras relate to
unitary representations of Clifford--Lie  {super}groups.

A finite dimensional (real) Lie superalgebra $\fn=\fn_0\op\fn_1$ is called a
{\it Clifford--Lie} superalgebra if
$\fn_0$ is a subset of the center of $\fn.$
Then $N := (\fn_0,+)$ is a Lie group with Lie algebra $\fn_0$ and
$\cN := (N,\fn)$ is a Lie supergroup for which $N$ acts trivially on $\fn$;
in \cite{CCTV} these  {super}groups are called {\it super translation groups}.

\begin{deft} {\rm We say that a unitary representation
$(\pi,\chi,\hH)$ of $\cN$ is {\it $\lambda$-admissible}
for $\lambda \in \fn_0^*$ if
\[ \pi(x) = e^{i\lambda(x)}\one \quad \mbox{ for } \quad x \in \fn_0.\] }
\end{deft}

If a $\lambda$-admissible unitary
representation exists, then Lemma~\ref{necessary}(ii) implies
that $\lambda \in \sc(\fn)^\star$. By Schur's Lemma, every irreducible
unitary representation of $\cN$ is $\lambda$-admissible for some
$\lambda \in \sc(\fn)^\star$.

Fix $\lam\in \sc(\fn)^\star$. Then
\begin{equation}
  \label{eq:mu-lam1}
\mu_{\lambda} \: \fn_1\times\fn_1\longrightarrow \bbbr,
\quad\quad \mu_\lambda(x,y) := \frac{1}{2}\lam([x,y])
\end{equation}
is a positive semidefinite symmetric bilinear form on $\fn_1$, hence
defines on the quotient space
\[ \fn_{1,\lambda} := \fn_1/\{ x \in \fn_1 \: \mu_\lambda(x,x) = 0\} \]
a positive definite form $\oline\mu_\lambda$.
We write
\begin{equation}
  \label{eq:cliff-lam}
 \cC_\lambda := \cC(\fn_{1,\lambda}, \oline\mu_{\lambda})
\end{equation}
for the corresponding Clifford algebra and $\cC_{\lambda,\C}$ for its complexification.
This is a $C^*$-algebra whose representations are precisely the complex linear
 extensions of selfadjoint representations of $\cC_\lambda$.
The corresponding structure map lifts to a linear map
\[ \iota \: \fn_1 \to \cC_\lambda, \quad x \mapsto  \oline x \]
satisfying
\[ \iota(x)^2 = \mu_\lambda(x,x)\one = \frac{1}{2}\lam([x,x])
\quad \mbox{ for } \quad x \in \fn_1.\]
Therefore the map
\begin{equation}
\label{iota}
\iota_\lam:\fn\longrightarrow \cC_{\lambda,\C}, \quad
x\mapsto
\begin{cases}
i\lam(x)\one & \text{ for } x\in \fn_0\\
e^{\frac{\pi i}{4}}\cdot x& \text{ for } x\in\fn_1
\end{cases}
\end{equation}
is a homomorphism of Lie superalgebras if $\cC_{\lambda,\C}$ in endowed with the
natural Lie superbracket defined on homogeneous elements by the super-commutator
\[ [a,b] := ab - (-1)^{  {|a|  |b|}} ba, \quad a,b \in \cC_{\lambda,\C}.\]
If $\rho$ is a selfadjoint representation of $\cC_{\lambda},$
then $\rho_\C\circ\iota_\lam$ is a $\lambda$-admissible
unitary representation of $\fn$, so that
$\pi := e^{i\lambda}$ on $\fn_0$ leads to a unitary
$\lambda$-admissible representation of the Lie supergroup $\cN$.
If, conversely, $(\pi, \chi_\pi)$ is  a
$\lambda$-admissible representation of $\cN$,
then the universal property of $\cC_\lambda$ implies the existence
of a selfadjoint representation
$\rho$ of $\cC_\lambda$ with $\rho_\C \circ \iota_\lambda = \chi_\pi$.
This establishes a  {one-to-one} correspondence between
$\lambda$-admissible unitary representations of $\cN$
and selfadjoint representations of $\cC_\lambda$, which
are completely described in Proposition~\ref{unique}.
We thus obtain:

\begin{Pro}\label{com-red}
{\rm($\lambda$-admissible representations)}  Suppose that $\lambda \in \sc(\fn)^\star$.
Then there exists a finite dimensional irreducible unitary
representation $(\pi_\lambda,\chi_\lambda, \hN_\lambda)$ of $\cN = (N,\fn)$
which is unique if $\dim \fn_{1,\lambda}$ is odd and unique up to parity reversal if
$\dim \fn_{1,\lambda}$ is even. Any other $\lambda$-admissible unitary representation
$(\pi,\hH)$ is of the form $\pi = \one \otimes \pi_\lambda$, where
$\hH = \hM \otimes \hN$ is a tensor product of Hilbert superspaces.
If $\dim \fn_{1,\lambda}$ is odd, then the multiplicity space $\hM$ can
be chosen purely even.
\end{Pro}

\begin{cor} {\rm(Irreducible  {representations})}
Each irreducible unitary representation of $\cN = (N,\fn)$
is finite dimensional and equivalent to some
$(\pi_\lambda, \chi_\lambda)$ with $\lambda \in \sc(\fn)^\star$
or to $\pi_\lambda^\Pi$, where $\Pi$ denotes parity reversal, i.e.,
$\pi_\lambda^\Pi = \pi_\lambda$ but $\cH_j^\Pi := \cH_{1-j}$ for $j =0,1$.
If $\dim \fn_{1,\lambda}$ is odd, then both are equivalent.
\end{cor}

\subsection{Semidirect products}

Now we assume that $K$ is a compact Lie  group acting on $\fn$
by a group homomorphism  $\rho:K\longrightarrow \Aut(\fn)
 {\cong \Aut(N)}.$
The action $\rho$ induces the action $\rho^\star$ of $K$ on the
{\it dual cone}
\[ \sc(\fn)^\star=\{\lam\in\fn_0^*\mid (\forall x \in \fn_1)\
 \lam([x,x])\geq 0\}
\quad \mbox{ by } \quad
\rho^\star_k(\lam):=\lam\circ \rho^{-1}_k,\]
for $k\in K$ and $\lam\in\sc(\fn)^\star.$
We write $\cO_\lambda$ for the corresponding orbit of $\lambda$.

We now explain how the irreducible unitary representations of the Lie supergroup
$\cG := \cN \rtimes K = ( {N}\rtimes K, \fn)$
can be classified. This classification has been derived in
\cite{CCTV} (Theorems~4 and 5) by generalizing Mackey's Imprimitivity
to the super context and by using the corresponding technique of unitary
induction (Definition~\ref{def:ind-rep}).
We now give a precise formulation of this result.

Fix $\lambda \in \sc(\fn)^\star$ and consider its stabilizer group
$K_\lambda \subeq K$. In the Mackey context,
$\cG_\lambda := \cN \rtimes K_\lambda$ plays the role of the {\it little group}
from which we want to induce  {representations} to~$\cG$. Therefore we first have
to classify the unitary representations of $\cG_\lambda$
which are $\lambda$-admissible in the sense that
$\pi(x) = e^{i\lambda(x)}\one$ for $x \in \fn_0$.
This can be done with the tools developed in the
preceding two subsections.

As $K_\lambda$ preserves $\lambda$,
its action on $\fn_1$ factors through an orthogonal representation
\[ \rho_\lambda \: K_\lambda \to \OO(\fn_{1,\lambda}).\]
Note that $K_\lambda$ need not be connected, so that the range of
$\rho_\lambda$ need not be contained in the identity component
$\SO(\fn_{1,\lambda})$. This causes some trouble in the constructions
because it leads to graded unitary representations.
Therefore a key point in the construction in \cite{CCTV} is to consider the
following subgroup of $K_\lambda$:
\begin{equation}\label{K-circle} K_\lambda^\circ :=
\begin{cases}
\rho_\lambda^{-1}(\SO(\fn_{1,\lambda})) & \text{ if }
\rho_\lambda(K_\lambda) \not\subset \SO(\fn_{1,\lambda})\mbox{ and }
\dim(\fn_{1,\lambda}) \ \text{even} \\
K_\lambda & \text{otherwise}.
\end{cases}
\end{equation}
So either $K_\lambda^\circ$ equals $K_\lambda$ or it is a subgroup of index~$2$.
From the central extension $\hat\OO(\fn_{1,\lambda})$ of
$\OO(\fn_{1,\lambda})$ by $\{\pm1\}$
(Proposition~\ref{prop:ohat}), we obtain a central extension
\[ \hat K_\lambda= \rho_\lambda^*\hat\OO(\fn_{1,\lambda})
= \{ (k,g)  \in K_\lambda \times \hat\OO(\fn_{1,\lambda}) \:
\rho_\lambda(k) = \oline g\}.\]

Recall the Clifford algebra
$\cC_\lambda := \cC(\fn_{1,\lambda}, \oline\mu_{\lambda})$
from \eqref{eq:cliff-lam} and its irreducible
representation $(\tau, \hN_\lambda)$.
We then obtain with Proposition~\ref{prop:ohat} a graded unitary
representation $\kappa_\lambda$ of $\hat K_\lambda$ on $\hN$ by
\[ \kappa_\lambda \: \hat K_\lambda \to \Aut(\hN_\lambda), \quad
\kappa_\lambda(k,g) := \kappa(g)\]
satisfying
\[  \kappa_\lambda(k,g)\tau(\iota(x)) \kappa(k,g)^{-1} = \tau(\iota(\rho_\lambda(k)x))
\quad \mbox{ for } \quad
(k,g) \in \hat K_\lambda, x \in \fn_1.\]
For the corresponding unitary representation $\chi_\lambda$ of $\fn$
on $\hN_\lambda$, this leads to
\[  \kappa_\lambda(k,g)\chi_\lambda(x)\kappa_\lambda(k,g)^{-1} = \chi_\lambda(\rho_\lambda(k)x)
\quad \mbox{ for } \quad
(k,g) \in \hat K_\lambda, x \in \fn,\]
so that it combines with $\kappa_\lambda$ to a unitary representation
$(\hat\pi_\lambda, \hat\chi_\lambda)$ of $\hat\cG_\lambda = \cN \rtimes \hat K_\lambda$
defined by
\[ \hat\pi_\lambda(x,k)
= e^{i\lambda(x)} \kappa_\lambda(k), \qquad
\hat\chi_\lambda(x,y)
= \chi_\lambda(x) + d\kappa_\lambda(y), \quad
x \in \fn, y \in \fk_\lambda.\]

We call a unitary representation $(\pi,\hH)$ of
$\hat K_\lambda^\circ$ {\it odd} if
$\pi(-1) = -\one$. Note that $\pi(-1)$ is always a unitary involution
and that $\pi(-1) \in \{\pm \one\}$ holds if
$\pi$ is irreducible by Schur's Lemma.

 {Theorem~4 of \cite{CCTV} now} asserts the existence of a functor
$r \mapsto \theta^\lambda_r$ which assigns to an odd unitary representation
$(r,\hH)$ of $\hat K_\lambda^\circ$ a unitary representation of
the {\it little supergroup} $\cN \rtimes K_\lambda$ which is $\lambda$-admissible.
This establishes in particular a bijection of the equivalence classes
of irreducible odd unitary representation of $\hat K_\lambda^\circ$
on (purely even) Hilbert spaces
with the irreducible unitary representation of the Lie supergroup
$\cN \rtimes K_\lambda$.

More concretely, starting from   {an odd unitary representation $r$ of $\hat K_\lambda^\circ,$} we first construct the induced
representation $(\hat r, \hat \hH)$ (Remark~\ref{rem:super-ext})
which is a graded representation
of the full group $\hat K_\lambda$ on a Hilbert superspace.
Now
\[ \theta^\lambda_r(x,k) := e^{i\lambda(x)}
\hat r(k) \otimes \kappa_\lambda(k) \]
defines the corresponding unitary representation of $ {N} \rtimes K_\lambda$
(here we use that $r$ and $\kappa_\lambda$ are both odd and graded with
respect to the subgroup $\hat K_\lambda^\circ$) and the associated
representation of the Lie superalgebra is determined by
\[ \chi(x) = \tau_\lambda(\iota_\lambda(x)) \quad \mbox{ for } \quad x \in \fn_1.\]
As $\theta^\lambda_r$ is a unitary representation of the little supergroup
$\cN \rtimes K_\lambda$, unitary induction now leads to a unitary represntation
\[ \Theta^\lambda_r := \Ind_{\cN \rtimes K_\lambda}^{\cN \rtimes K}(\theta^\lambda_r) \]
of the Lie super $\cN \rtimes K$.

\begin{Thm}\label{main1} Fix $\lam\in \sc(\fn)^\star$.
If $(r,\hH)$ is an odd unitary representation of
$\hat K_\lambda^\circ$, then $\Theta^\lambda_r$
is a unitary representation of the Lie supergroup $\cN \rtimes K$.
If $r$ is irreducible, then so is $\Theta^\lambda_r$ and
all irreducible unitary representations of $\cN \rtimes K$ are obtained  {in?}
this way. If $\Theta^\lambda_r$ and $\Theta^{\lambda'}_{r'}$ are equivalent,
then $\lambda'$ and $\lambda$ lie in the same $K$-orbit in $\sc(\fn)^\star$
and $\Theta^\lambda_r \cong \Theta^{\lambda}_{r'}$
if and only if $r \cong r'$. Hence the equivalence classes of irreducible
unitary representations of $\cN \rtimes K$
are determined by the $K$-orbits $\cO_\lambda$ in $\sc(\fn)^\star$
and, for any fixed $\lambda$, they are in {one-to-one} correspondence
with odd irreducible representations of $\hat K_\lambda^\circ$.
\end{Thm}

\end{document}